\documentclass[a4paper,twoside,10pt,reqno,final]{amsart}
\usepackage{amssymb}
\usepackage[cm,headings]{fullpage}
\usepackage[kerning=true,spacing=true,final]{microtype}

\usepackage[utf8]{inputenc}
\usepackage[T1]{fontenc}

\usepackage{enumerate}

\let\comp=\circ
\let\wtilde=\widetilde

\usepackage[final]{hyperref}

\newtheorem{theorem}{Theorem}
\newtheorem{corollary}[theorem]{Corollary}
\newtheorem{proposition}[theorem]{Proposition}
\newtheorem{lemma}[theorem]{Lemma}
\newtheorem{fact}[theorem]{Fact}
\newtheorem*{th_ex_om}{Theorem~\ref{t:ext_C1om}}
\theoremstyle{definition}
\newtheorem{definition}[theorem]{Definition}
\theoremstyle{remark}
\newtheorem{remark}[theorem]{Remark}
\newtheorem{example}[theorem]{Example}

\newcommand\setsep{;\allowbreak\ } 

\newcommand\abs[1]{\mathopen|#1\mathclose|}
\newcommand\absa[1]{\left|#1\right|}
\newcommand\absb[2]{\csname#1l\endcsname|#2\csname#1r\endcsname|}
\newcommand\norm[1]{\mathopen\|#1\mathclose\|}

\newcommand\normb[2]{\csname#1l\endcsname\|#2\csname#1r\endcsname\|}
\newcommand\tnorm[1]{\mathopen{|\mskip-0.5\thinmuskip|\mskip-0.5\thinmuskip|}#1\mathclose{|\mskip-0.5\thinmuskip|\mskip-0.5\thinmuskip|}}

\newcommand{\cl}[2][3]{{}\mkern#1mu\overline{\mkern-#1mu#2}}

\newcommand\al{\alpha}
\newcommand\ga{\gamma}
\newcommand\Ga{\Gamma}
\newcommand\de{\delta}
\newcommand\ve{\varepsilon}
\newcommand\vp{\varphi}
\newcommand\om{\omega}

\newcommand\N{{\mathbb N}}
\newcommand\R{{\mathbb R}}

\renewcommand\d{\,\mathrm{d}}
\newcommand\restr[1]{\mathclose\restriction_{#1}}

\newcommand\mc{\mathcal}

\newcommand\ev[2]{\langle #1,#2\rangle}

\newcommand\evb[3]{\csname#1l\endcsname\langle #2,#3\csname#1r\endcsname\rangle}

\DeclareMathOperator{\supp}{supp}
\DeclareMathOperator{\suppo}{supp_o}
\DeclareMathOperator{\dist}{dist}
\DeclareMathOperator{\diam}{diam}
\DeclareMathOperator{\len}{len}

\newcommand\Cloc[1]{C^{#1}_{\mathrm{loc}}}
\newcommand\Cob[1]{C^{#1}_{\mathrm B}}

\newcommand\eqdef{\mathrel{\mathop:}=}

\newcommand\cWG[1]{(\mathrm{WG}_{#1})}
\newcommand\cWd[1]{(\mathrm W_{#1})}
\newcommand\cW{(\widetilde{\mathrm W})}

\begin{document}
\title[On Whitney-type extension theorems on Banach spaces]{On Whitney-type extension theorems on Banach spaces for $C^{1,\om}$, $C^{1,+}$, $\Cloc{1,+}$, and $\Cob{1,+}$-smooth functions}
\author{Michal Johanis}
\address{Charles University, Faculty of Mathematics and Physics\\Department of Mathematical Analysis\\Sokolovská~83\\186~75 Praha~8\\Czech Republic}
\email{johanis@karlin.mff.cuni.cz}
\author{Václav Kryštof}
\address{Charles University, Faculty of Mathematics and Physics\\Department of Mathematical Analysis\\Sokolovská~83\\186~75 Praha~8\\Czech Republic}
\author{Luděk Zajíček}
\address{Charles University, Faculty of Mathematics and Physics\\Department of Mathematical Analysis\\Sokolovská~83\\186~75 Praha~8\\Czech Republic}
\email{zajicek@karlin.mff.cuni.cz}
\date{March 2023}
\begin{abstract}
Our paper is a complement to a recent article by D.~Azagra and C.~Mudarra (2021).
We show how older results on semiconvex functions with modulus~$\om$ easily imply extension theorems for $C^{1,\om}$-smooth functions on super-reflexive Banach spaces which are versions of some theorems of Azagra and Mudarra.
We present also some new interesting consequences which are not mentioned in their article, in particular extensions of $C^{1,\om}$-smooth functions from open quasiconvex sets.
They proved also an extension theorem for $\Cob{1,+}$-smooth functions (i.e., functions with uniformly continuous derivative on each bounded set) on Hilbert spaces.
Our version of this theorem and new extension results for $C^{1,+}$ and $\Cloc{1,+}$-smooth functions (i.e., functions with uniformly, resp. locally uniformly continuous derivative),
all of which are proved on arbitrary super-reflexive Banach spaces, are further main contributions of our paper.
Some of our proofs use main ideas of the article by D.~Azagra and C.~Mudarra, but all are formally completely independent on their article.
\end{abstract}
\maketitle

\section{Introduction}\label{sec:intro}

Our paper is a complement to a recent article~\cite{AM} which concerns extension theorems for $C^{1,\om}$ and $\Cob{1,+}$-smooth functions (see Definitions~\ref{d:C1om} and~\ref{d:C1+}) on Hilbert and super-reflexive Banach spaces.

Our main contributions are extension theorems for $C^{1,+}$ (see Definition~\ref{d:C1+}) and $\Cob{1,+}$-smooth functions (Theorems~\ref{t:ext_C1+_gen}, \ref{t:ext_lipC1+} and~\ref{t:ext_CB1+}, resp.~\ref{t:ext_C1ob}) which in particular generalise results from~\cite{AM} (see~\cite[Theorem~6.1 and its proof]{AM})
proved in a Hilbert space to an arbitrary super-reflexive Banach space
and a new extension theorem for $\Cloc{1,+}$-smooth functions (Theorem~\ref{t:ext_C1+loc}) which is also proved in an arbitrary super-reflexive Banach space.

Our contribution concerning the $C^{1,\om}$-smooth extensions from arbitrary sets is not big, since our theorems relatively easily follow from corresponding results of~\cite{AM}, see Remark~\ref{r:relation_to_AM1}(d).
However, we believe that our proofs of Theorems~\ref{t:ext_C1om} and~\ref{t:extC1_Holder} which are ``qualitative'' versions of several theorems from~\cite{AM}
(which are ``quantitative'', see Remark~\ref{r:relation_to_AM1}(d)) are worth publishing since they are very short and transparent
as they use older results on $\om$-semiconvex functions from~\cite{DZ1} and~\cite{Kr1} (see Lemma~\ref{l:smooth_nu} and Theorem~\ref{t:Ilmanen} below).
Note that these results on $\om$-semiconvex functions are independently (implicitly) proved in~\cite{AM} (in a slightly different setting) using a notion of a ``strongly $C\vp$-paraconvex function'',
which does not correspond to Rollewicz's notion of strong paraconvexity (cf. Remark~\ref{r:paraconvex}) but coincides with $\al$-semiconvexity, where $\al(t)=C\vp(t)/t$ for $t>0$ and $\al(0)=0$.
Nevertheless, the main idea of the proofs of Theorems~\ref{t:ext_C1om} and~\ref{t:extC1_Holder} is taken from~\cite{AM}, namely we use the crucial observation from \cite{AM} (cf. Lemma~\ref{l:WG->envelopes} and the remark before it)
that natural necessary conditions for an extension imply the inequality $h\le H$ from Lemma~\ref{l:WG->envelopes} (or $m\le g$ from~\cite[p.~9]{AM}).

We use the classical Whitney-Glaeser condition $\cWG\om$ (see Theorem~\ref{t:ext_C1om}) in contrast to conditions $A(f,G)<+\infty$, $(W^{1,\om})$ and $(mg^{1,\om})$ used in~\cite{AM}.
Note that the use of condition $\cWG\om$ enables us to naturally apply Whitney-type extension results to ``true'' extension theorems for $C^{1,\om}$-smooth functions not only on open convex sets,
but also on open \emph{quasiconvex} sets (Corollaries~\ref{c:extC1om-quasic}, \ref{c:extC1_Hold-quasic}).

We will use the following basic notions:
\begin{definition}\label{d:C1om}
\hfill
\begin{enumerate}[(i)]
\item
We denote by $\mc M$ the set of all finite moduli, i.e. the functions $\om\colon[0,+\infty)\to[0,+\infty)$ which are non-decreasing and continuous at $0$ with $\om(0)=0$.

\item
If $U$ is an open subset of a normed linear space and $\om\in\mc M$,
we denote by $C^{1,\om}(U)$ the set of all Fréchet differentiable $f\colon U\to\R$ such that the Fréchet derivative $Df$ is uniformly continuous with modulus $C\om$ for some $C\ge0$,
that is
\[
\norm{Df(x)-Df(y)}\le C\om(\norm{x-y})\quad\text{whenever $x,y\in U$.}
\]
\end{enumerate}
\end{definition}
Using results of~\cite{DZ1} and~\cite{Kr1} we easily obtain the following result which can be deduced from~\cite{AM} (see Remark~\ref{r:relation_to_AM1}(d)),
where, however, no equivalent version is explicitly formulated.

\begin{theorem}\label{t:ext_C1om}
Let $\om\in\mc M$ be a concave modulus and let $X$ be a super-reflexive Banach space that has an equivalent norm with modulus of smoothness of power type $2$.
Let $E\subset X$ and let $f$ be a real function on $E$.
Then $f$ can be extended to a function $F\in C^{1,\om}(X)$ if and only if the following condition $\cWG\om$ holds:

There exist a mapping $G\colon E\to X^*$ and $M>0$ such that
\begin{align*}
\norm{G(y)-G(x)}&\le M\om(\norm{y-x}),\\
\absb{big}{f(y)-f(x)-\ev{G(x)}{y-x}}&\le M\om(\norm{y-x})\norm{y-x}
\end{align*}
for each $x,y\in E$.

Moreover, if $\cWG\om$ is satisfied, then $F$ can be found such that $DF(x)=G(x)$ for each $x\in E$.
\end{theorem}
The assumption on $X$ cannot be relaxed, see Remark~\ref{r:part_unity}(a).
By Remark~\ref{r:Lp-modulus}, Theorem~\ref{t:ext_C1om} (and its Corollary~\ref{c:extC1om-quasic}) can be applied e.g. if $X$ is isomorphic to some $L_p(\mu)$ with $p\ge2$.

\begin{remark}\label{r:relation_to_AM1}
\hfill

(a)
Condition $\cWG\om$ naturally corresponds to Whitney's condition in his celebrated theorem from \cite{Wh1} on extensions to $C^k$-smooth functions in $X=\R^n$ (for $k=1$)
and appears explicitly in Glaeser's thesis~\cite{Gl}.

(b)
It was stated several times that Theorem~\ref{t:ext_C1om} was proved for $X=\R^n$ by Glaeser in~\cite{Gl} (see e.g. \cite[p.~516]{BS}, \cite[p.~2]{AM}).
However,~\cite{Gl} (which works with $C^{k,\om}$-smooth functions) contains only a corresponding result on extension of a function $f$ defined on a closed subset $F$ of a compact interval $K\subset\R^n$ (``pavé compact'') to a function belonging to $C^{1,\om}(K)$
(see~\cite[Proposition VII, p.~33, $m=1$]{Gl}; cf. corresponding \cite[Complement~3.6, p.~9]{Ma}).

(c)
Whitney-Glaeser condition $\cWG\om$ is equivalent to condition~\cite[(1.2)]{AM} (for $X=\R^n$), which we obtain from $\cWG\om$ by changing ``$\le M\om(\norm{y-x})\norm{y-x}$'' to ``$\le M\vp(\norm{y-x})$'',
where $\vp(t)=\int_0^t\om(s)\d s$ for $t\ge0$ (cf.~Fact~\ref{f:om_vp} below).

(d)
If $X$ is a Hilbert space, then the statement of Theorem~\ref{t:ext_C1om} follows from~\cite[Theorem 1.6]{AM}.
Moreover, in full generality it follows from~\cite[Theorem 1.9]{AM} using our Lemma~\ref{l:WG->envelopes}, Lemma~\ref{l:equiv}, and the easy observation from Remark~\ref{r:relation_to_AM2} below.
Further, results in~\cite{AM} assert not only that $C^{1,\om}$-smooth extension $F$ exists, but they are quantitative in the sense that
they assert that $DF$ is uniformly continuous with modulus $K\om$, where $K>0$ (which does not depend on $f$) is explicitly given.
Examining the proofs of Theorem~\ref{t:Ilmanen} and Lemma~\ref{l:smooth_nu} we could obtain similar quantitative results for Theorems~\ref{t:ext_C1om} and~\ref{t:extC1_Holder} too, but less precise than those in~\cite{AM} which seem to be close to the (unknown) optimal ones.
\end{remark}

We moreover prove Theorem~\ref{t:extC1_Holder}, an extension theorem for $C^{1,\al}$-smooth functions which works e.g. in $L_p(\mu)$ spaces with $p\ge1+\al$.
Theorem~\ref{t:extC1_Holder} can be also deduced from \cite{AM}.
However, our (simple but important) observation that the assumption on $X$ in Theorem~\ref{t:ext_C1om} and Theorem~\ref{t:extC1_Holder} cannot be relaxed (see Remark~\ref{r:part_unity}) is new.

We are also interested in the extension of $C^{1,\om}$-smooth functions from open sets, which is not considered in~\cite{AM}.
As a consequence of Theorems~\ref{t:ext_C1om} and~\ref{t:extC1_Holder} we obtain new interesting extension theorems from so-called open quasiconvex (and in particular convex) sets.
These results do not use condition $\cWG\om$ and so they are true extension theorems for $C^{1,\om}$-smooth functions.
(The set is called quasiconvex if its inner metric is bi-Lipschitz equivalent to the metric of the enclosing space; see Definition~\ref{d:quasic}.
Quasiconvex sets are sometimes called ``sets with Whitney arc property'', since they were used in Whitney's article~\cite{Wh2} concerning extensions of $C^k$-smooth functions in $\R^n$.)
Namely, in Corollary~\ref{c:extC1om-quasic} we prove that if $X$ is as in Theorem~\ref{t:ext_C1om} and $\om\in\mc M$,
then each function from $C^{1,\om}(U)$, $U\subset X$ quasiconvex, can be extended to a function from $C^{1,\om}(X)$; see also analogous Corollary~\ref{c:extC1_Hold-quasic}.

Another new result is Theorem~\ref{t:ext_C1ob} on extensions for functions from $\Cob{1,+}(X)$ (i.e. functions with uniformly continuous derivative on every bounded set),
which is in fact a generalisation of~\cite[Theorem 6.1]{AM} from Hilbert spaces to arbitrary super-reflexive spaces.
A simplified formulation is as follows:
\begin{theorem}\label{t:ext_CB1+}
Let $X$ be a super-reflexive Banach space.
Let $E\subset X$ and let $f$ be a real function on $E$.
Then $f$ can be extended to a real function $F\in\Cob{1,+}(X)$ if and only if the following condition $\cW$ holds:

There exists a mapping $G\colon E\to X^*$ such that for each bounded $B\subset E$ the function $f$ is bounded on $B$, the mapping $G$ is bounded and uniformly continuous on $B$,
and for each $\ve>0$ there exists $\de>0$ such that $\absb{big}{f(y)-f(x)-\ev{G(x)}{y-x}}\le\ve\norm{y-x}$ whenever $x,y\in B$ and $\norm{y-x}<\de$.

Moreover, if $\cW$ is satisfied, then $F$ can be found such that $DF(x)=G(x)$ for each $x\in E$.
\end{theorem}

We also prove two results (Theorems~\ref{t:ext_C1+_gen} and ~\ref{t:ext_lipC1+}) on extensions from arbitrary sets for $C^{1,+}$-smooth functions.
As a corollary we obtain extension results for $C^{1,+}$-smooth functions on open quasiconvex sets (Corollaries~\ref{c:extLipC1+-quasic}, \ref{c:extC1+-quasic}, and~\ref{c:extC1+B-quasic}).
Another interesting theorem that has no analogue in~\cite{AM} is Theorem~\ref{t:ext_C1+loc},
which in particular characterises those functions on a closed subset of a super-reflexive space that can be extended to functions on the whole space having locally uniformly continuous derivative.

The structure of the present article is following:
In (unfortunately rather long) Section~\ref{sec:prelim} we fix the notation and recall the needed known facts on moduli, $\om$-semiconvex functions, and super-reflexive spaces.
The core of the article is contained in a short Section~\ref{sec:basic}, where we prove two basic lemmata and one proposition.
Using results of Section~\ref{sec:basic} we rather easily prove Theorem~\ref{t:ext_C1om} and other results concerning the $C^{1,\om}$ case in Section~\ref{sec:C1om}.
In Sections~\ref{sec:C1+} and~\ref{sec:CB1+} we prove new results concerning $C^{1,+}$, $\Cob{1,+}$, $\Cloc{1,+}$ cases, in particular Theorem~\ref{t:ext_C1ob} (which contains Theorem~\ref{t:ext_CB1+}).

We remark that all the positive extension results are proved in super-reflexive spaces.
Remarks~\ref{r:part_unity}, \ref{r:quasic-bump} and \ref{r:ext_loc} show that for all these results this assumption on~$X$ is necessary.

Some natural questions remain open, see Remarks~\ref{r:moduli_nonconc}, \ref{r:ext-gen}, \ref{r:quasic_unbound}, and \ref{r:ext_convex}.

\section{Preliminaries}\label{sec:prelim}

\subsection{Basic notation, properties of moduli, and quasiconvex sets}

All the normed linear spaces considered are real.
Let $(X,\norm\cdot)$ be a normed linear space.
By $U(x,r)\subset X$, resp. $S(x,r)\subset X$ we denote the open ball, resp. sphere centred at $x\in X$ with radius $r>0$.
By $B_X$, resp. $S_X$ we denote the closed unit ball, resp. unit sphere in~$X$.
For $x\in X$ and $g\in X^*$ we will denote the evaluation of $g$ at $x$ also by $\ev gx=g(x)$.
An $L$-Lipschitz mapping is a mapping with a (not necessarily minimal) Lipschitz constant~$L$.
By $Df(x)$ we will denote the Fréchet derivative of $f$ at $x$, its evaluation in the direction $h$ will be denoted by $Df(x)[h]$.
We remind that the set of all finite moduli $\mc M$ and the $C^{1,\om}$-smooth functions are defined in Definition~\ref{d:C1om}.
Sometimes we will work with several norms on a space.
If it is necessary to specify objects with respect to a particular norm we will use the notation $U_{\norm\cdot}(0,r)$, $C_{\norm\cdot}^{1,\om}$, etc.

For a mapping $f\colon X\to Y$, where $X$ is a set and $Y$ is a vector space, we denote $\suppo f=f^{-1}(Y\setminus\{0\})$.
If $X$ is a topological space, then we denote $\supp f=\overline{\suppo f}$.

For the (standard) definition and properties of the length (variation) of a continuous curve $\ga\colon [0,1]\to X$ in a metric space $X$ see e.g.~\cite{Ch}.

Recall that a system $\{\psi_\al\}_{\al\in\Lambda}$ of functions on a set $X$ is called a partition of unity if
\begin{itemize}
\item $\psi_\al\colon X\to[0,1]$ for all $\al\in\Lambda$,
\item $\sum\limits_{\al\in\Lambda}\psi_\al(x)=1$ for each $x\in X$.
\end{itemize}
We say that the partition of unity $\{\psi_\al\}_{\al\in\Lambda}$ is subordinated to a covering $\mc U$ of $X$ if $\{\suppo\psi_\al\}_{\al\in\Lambda}$ refines $\mc U$,
i.e. for each $\al\in\Lambda$ there is $U\in\mc U$ such that $\suppo\psi_\al\subset U$.
Further, in case that~$X$ is a topological space we say that the partition of unity $\{\psi_\al\}_{\al\in\Lambda}$ is locally finite if the system $\{\suppo\psi_\al\}_{\al\in\Lambda}$ is locally finite,
i.e. if for each point $x\in X$ there is a neighbourhood $U$ of $x$ such that the set $\{\al\in\Lambda\setsep\suppo\psi_\al\cap U\neq\emptyset\}$ is finite.

\begin{definition}
Let $(X,\rho)$, $(Y,\sigma)$ be metric spaces and $f\colon X\to Y$.
Then:
\begin{enumerate}[(i)]
\item
The \emph{minimal modulus of continuity of $f$} is the function $\om_f\colon[0,+\infty)\to[0,+\infty]$ defined in the usual way:
\[
\om_f(\de)=\sup\bigl\{\sigma(f(x),f(y))\setsep x,y\in X,\rho(x,y)\le\de\bigr\}.
\]
\item
We say that $f$ is \emph{uniformly continuous with modulus $\om\in\mc M$} if $\sigma(f(x),f(y))\le\om(\rho(x,y))$ for every $x,y\in X$, i.e. $\om_f\le\om$.
\item
We say that $f$ is \emph{$\al$-Hölder} ($\al>0$) if $f$ is uniformly continuous with modulus $\om(t)=Ct^\al$ for some $C>0$.
\end{enumerate}
\end{definition}

Clearly, $\om_f$ is non-decreasing, and the mapping $f$ is uniformly continuous if and only if $\lim_{t\to0+}\om_f(t)=0$.

Recall that $\om\colon[0,+\infty)\to[0,+\infty]$ is \emph{sub-additive} if $\om(s+t)\le\om(s)+\om(t)$ for any $s,t\in[0,+\infty)$.
For $\om$ sub-additive it immediately follows by induction that
\begin{equation}\label{e:subad_Nt}
\om(Nt)\le N\om(t)\quad\text{for each $t\ge0$, $N\in\N$.}
\end{equation}
Consequently, the following fact holds:
\begin{fact}\label{f:subad-finite}
If $\om\colon[0,+\infty)\to[0,+\infty]$ is sub-additive and finite on a neighbourhood of $0$, then it is finite everywhere.
\end{fact}

The following fact is easy and well-known (it follows from~\cite[Remark~1, p.~407]{AP} and Fact~\ref{f:subad-finite}):
\begin{fact}\label{f:conv_subad}
If $A$ is a convex subset of a normed linear space, $Y$ is a metric space, and $f\colon A\to Y$ is uniformly continuous, then $\om_f$ is sub-additive and $\om_f\in\mc M$.
\end{fact}

The next two standard facts are also very easy to see.
\begin{fact}\label{f:concave-subad}
If $\om\in\mc M$ is concave, then it is sub-additive and $\om(Kt)\le K\om(t)$ for any $t\ge0$, $K\ge1$.
\end{fact}

\begin{fact}\label{f:common_mod}
If $C$ is a convex subset of a normed linear space, $Y$ is a metric space, and $f_1,\dotsc,f_n\colon C\to Y$ are uniformly continuous,
then there is $\om\in\mc M$ such that all $f_1,\dotsc,f_n$ are uniformly continuous with modulus~$\om$.
\end{fact}
\begin{proof}
It suffices to set $\om=\om_{f_1}+\dotsb+\om_{f_n}$ and notice that $\om\in\mc M$ by Fact~\ref{f:conv_subad}.
\end{proof}

We will also use the following well-known facts.

\begin{lemma}\label{l:Steckin}
Let $\om\in\mc M$.
\begin{enumerate}[(i)]
\item
If $\om$ satisfies $\limsup_{t\to+\infty}\frac{\om(t)}t<+\infty$, then there exists a concave $\wtilde\om\in\mc M$ such that $\om\le\wtilde\om$.
\item
If $\om$ is bounded, then there exists a bounded concave $\wtilde\om\in\mc M$ such that $\om\le\wtilde\om$.
\item
If $\om$ is sub-additive, then there exists a concave $\wtilde\om\in\mc M$ such that $\om\le\wtilde\om\le2\om$.
\end{enumerate}
\end{lemma}
The proof of (i) is quite easy, see~\cite[proof of Theorem~1, p.~407]{AP}, where the theorem states sub-additivity of $\wtilde\om$ but in fact the proof shows even concavity.
Alternatively one can use the statement of \cite[Theorem~1, p.~406]{AP} together with (iii).
Statement (ii) follows easily from (i).
Fact (iii) is due to Stechkin, but was not published by him.
For a proof, see~\cite[p.~78]{E} or~\cite[p.~670]{Ko}.

\begin{fact}\label{f:om_vp}
Let $\om\in\mc M$.
If we set $\vp(t)=\int_0^t\om(s)\d s$ for $t\ge0$, then $\vp\in\mc M$ and $\vp$ is convex.
If moreover $\om$ is concave, then $t\om(t)\le2\vp(t)$ for every $t\ge0$.
\end{fact}
\begin{proof}
The convexity of $\vp$ follows from the fact that $\om$ is non-decreasing, see e.g.~\cite[p.~11, p.~13]{RV}.
Now suppose that $\om$ is concave.
Then $\frac{\om(t)}t\le\frac{\om(s)}s$ for each $0<s\le t$, and so
\[
t\om(t)=2\int_0^t\frac st\om(t)\d s\le2\int_0^t\om(s)\d s=2\vp(t).
\]
\end{proof}

\begin{definition}
A metric space $(X,\rho)$ is called \emph{almost convex}
if for any $x,y\in X$ and $s,t>0$ such that $\rho(x,y)<s+t$ there is $z\in X$ such that $\rho(x,z)\le s$ and $\rho(z,y)\le t$.
\end{definition}

Note that this notion goes back to Aronszajn's (unpublished) thesis; see \cite[p. 417]{AP}, where the name ``almost 3-hyperconvex metric space'' is used.
The following (almost obvious) fact is mentioned in \cite[p.~438, footnote~20]{AP}.

\begin{fact}\label{f:in_metric-alm_convex}
Let $(X,\rho)$ be a metric space where every two points can be connected by a rectifiable curve and let $\rho_{\rm i}$ be the corresponding inner (intrinsic) metric.
Then $(X,\rho_{\rm i})$ is an almost convex metric space.
\end{fact}

The article \cite{BV} contains (without a proof) the following interesting and probably new observation which we will use below in the proof of Lemma~\ref{l:quasic-concave_modulus}.
(Note that the assumption of uniform continuity of $f$ was forgotten in~\cite{BV}.
We have been informed by Taras Banakh that this assumption cannot be dropped, cf. Remark~\ref{r:alm_convex-uniform_cont}.)

\begin{proposition}\label{p:alm_convex-mod_subad}
Let $(X,\rho)$, $(Y,\sigma)$ be metric spaces and let $f\colon X\to Y$ be uniformly continuous.
If $X$ is almost convex, then $\om_f$ is sub-additive.
\end{proposition}
\begin{proof}
Let $a,b\in[0,+\infty)$ and $\ve>0$.
Let $\de>0$ be such that $\sigma(f(u),f(v))\le\ve$ whenever $\rho(u,v)\le\de$.
Pick any $x,y\in X$ such that $\rho(x,y)\le a+b<a+b+\frac\de2$.
Using the almost convexity we can find $u\in X$ such that $\rho(x,u)\le a$ and $\rho(u,y)\le b+\frac\de2<b+\de$.
Using the almost convexity again we can find $v\in X$ such that $\rho(u,v)\le\de$ and $\rho(v,y)\le b$.
Hence
\[
\sigma(f(x),f(y))\le\sigma(f(x),f(u))+\sigma(f(u),f(v))+\sigma(f(v),f(y))\le\om_f(a)+\ve+\om_f(b).
\]
It follows that $\om_f(a+b)\le\om_f(a)+\om_f(b)+\ve$, and since this holds for any $\ve>0$, the sub-additivity follows.
\end{proof}

\begin{remark}\label{r:alm_convex-uniform_cont}
If $(X,\rho)$ is totally convex (which is a stronger property then almost convexity), then the assumption of uniform continuity of $f$ is unnecessary, see \cite[Remark~1]{AP}.

However, the assumption of uniform continuity of $f$ cannot be dropped not only in Proposition~\ref{p:alm_convex-mod_subad} but also in the case when $(X,\rho)$ is an inner metric space.
Indeed let $X=\R^2\setminus\{(0,0)\}$ be equipped with the euclidean metric $\rho$ which clearly coincides with $\rho_{\rm i}$.
Set $f(-1,0)=0$; $f(x,y)= 1$ if $x\le0$, $(x,y)\neq(-1,0)$, $(x,y)\neq(0,0)$; $f(x,y)=2$ if $x>0$, $(x,y)\neq(1,0)$; and $f(1,0)=3$.
Then clearly $\om_f(2)=3>2=\om_f(1)+\om_f(1)$.
\end{remark}

The following notion of a quasiconvex space (or domain in $\R^n$) is nowadays a standard tool in Geometric Analysis.

\begin{definition}\label{d:quasic}
We say that a metric space $(X,\rho)$ is \emph{$c$-quasiconvex} for $c\ge1$ if for each $x,y\in X$ there exists a continuous rectifiable curve $\ga\colon[0,1]\to X$ such that
$\ga(0)=x$, $\ga(1)=y$, and $\len\ga\le c\rho(x,y)$, where $\len\ga$ is the length (variation) of the curve~$\ga$.
We say that $X$ is \emph{quasiconvex} if it is $c$-quasiconvex for some $c\ge1$.
\end{definition}

Note that convex subsets of normed linear spaces are $1$-quasiconvex.

\begin{remark}\label{r:biLip-quasic}
It is well-known and easy to prove that each bi-Lipschitz image of a quasiconvex metric space is quasiconvex.
\end{remark}

The following example, which is essentially well-known (cf. \cite[Lemma~5.9]{Va1}), will be used later in Remark~\ref{r:quasic-bump}.

\begin{example}\label{ex:shell-quasic}
Let $X$ be a normed linear space with $\dim X\ge2$ and $0\le r<R\le+\infty$.
Then the spherical shell $D=\{x\in X\setsep r<\norm x<R\}$ is $2$-quasiconvex.
Indeed, let $x,y\in D$.
We may assume that $\norm y\ge\norm x$.
Let $t=\max\{s\in[0,1]\setsep\norm{x+s(y-x)}\le\norm x\}$.
By \cite[Theorem~4J]{Sch} there is a continuous rectifiable curve $\ga$ joining $x$ and $x+t(y-x)$ in $S(0,\norm x)\subset D$ such that $\len\ga\le2t\norm{y-x}$.
By concatenating $\ga$ with the segment $[x+t(y-x),y]\subset D$ we obtain a continuous rectifiable curve joining $x$ and $y$ in $D$ whose length is at most $2t\norm{x-y}+(1-t)\norm{x-y}\le2\norm{x-y}$.
\end{example}

There is a plenty of non-convex quasiconvex sets, e.g. bi-Lipschitz images of convex sets and spherical shells are quasiconvex (by Remark~\ref{r:biLip-quasic}).

\begin{lemma}\label{l:quasic-concave_modulus}
Let $A$ be a $c$-quasiconvex subset of a normed linear space and let $f\colon A\to Y$ be a uniformly continuous mapping into a metric space $(Y,\sigma)$.
Then there exists a concave $\om\in\mc M$ such that $\frac12\om\le\om_f\le c\om$.
\end{lemma}
\begin{proof}
Let $\rho_{\rm i}$ be the inner metric on $A$ induced by the norm.
Note that $\norm{x-y}\le\rho_{\rm i}(x,y)\le c\norm{x-y}$ for each $x,y\in A$.
Then the mapping $f\colon(A,\rho_{\rm i})\to Y$ is uniformly continuous; denote by $\om^{\rm i}_f$ its minimal modulus of continuity.
Proposition~\ref{p:alm_convex-mod_subad} together with Fact~\ref{f:in_metric-alm_convex}, and Fact~\ref{f:subad-finite} imply that $\om^{\rm i}_f$ is sub-additive and finite.
Now Lemma~\ref{l:Steckin}(iii) implies that there is a concave $\om\in\mc M$ such that $\om^{\rm i}_f\le\om\le2\om^{\rm i}_f$.
Since $\norm{x-y}\le\rho_{\rm i}(x,y)$, it follows that $\frac12\om\le\om^{\rm i}_f\le\om_f$.
Further, using Fact~\ref{f:concave-subad} we obtain
\[
\sigma(f(x),f(y))\le\om^{\rm i}_f(\rho_{\rm i}(x,y))\le\om(\rho_{\rm i}(x,y))\le\om(c\norm{y-x})\le c\om(\norm{y-x})
\]
for each $x,y\in A$ and so $\om_f\le c\om$.
\end{proof}

The following application of quasiconvexity is well-known:
\begin{proposition}\label{p:om_der-lip-quasic}
Let $X$, $Y$ be normed linear spaces, $U\subset X$ an open $c$-quasiconvex set, and let $f\colon U\to Y$ be a Fréchet differentiable mapping.
If $K=\sup_{x\in U}\norm{Df(x)}<+\infty$, then $f$ is $cK$-Lipschitz.
\end{proposition}
\begin{proof}
The convexity of neighbourhoods of points in $U$ implies that $f$ is locally $K$-Lipschitz.
The statement now follows from \cite[Lemma~5.5]{Va2}.
\end{proof}

\subsection{\texorpdfstring{Subclasses of $C^1$-smooth functions and Whitney-type conditions}{Subclasses of C1-smooth functions and Whitney-type conditions}}

Besides the class of $C^{1,\om}$-smooth functions we will consider also the following subclasses of $C^1$-smooth functions:

\begin{definition}\label{d:C1+}
Let $X$ be a normed linear space and $U\subset X$ open.
We denote by $C^{1,\al}(U)$, $0<\al\le1$, the set of all Fréchet differentiable functions $f$ on $U$ such that $Df$ is $\al$-Hölder on $U$, i.e. $C^{1,\al}(U)=C^{1,\om}(U)$ for $\om(t)=t^\al$.
We denote by $C^{1,+}(U)$ the set of all Fréchet differentiable functions $f$ on $U$ such that $Df$ is uniformly continuous on $U$,
and by $\Cob{1,+}(U)$ the set of all Fréchet differentiable functions $f$ on $U$ such that $Df$ is uniformly continuous on each bounded subset of $U$.
Further, by $\Cloc{1,+}(U)$ we denote the set of all functions on $U$ which are locally $C^{1,+}$-smooth on $U$.
\end{definition}

For the reader's convenience we now recall two Whitney-type conditions $\cWG\om$ and $\cW$ used in Theorem~\ref{t:ext_C1om} and~\ref{t:ext_CB1+} in Section~\ref{sec:intro}
and introduce a third one, $\cWd G$.
\begin{definition}\label{d:condW}
Let $X$ be a normed linear space, $E\subset X$, $\om\in\mc M$, and $f\colon E\to\R$.
We say that $f$ satisfies condition
\smallskip
\begin{itemize}
\item[$\cWd G$]
if $G\colon E\to X^*$ is uniformly continuous and for each $\ve>0$ there exists $\de>0$ such that
\[
\absb{big}{f(y)-f(x)-\ev{G(x)}{y-x}}\le\ve\norm{y-x}
\]
whenever $x,y\in E$ and $\norm{y-x}<\de$;

\smallskip
\item[$\cWG\om$]
if there exist a mapping $G\colon E\to X^*$ and $M>0$ such that
\begin{align*}
\norm{G(y)-G(x)}&\le M\om(\norm{y-x}),\\
\absb{big}{f(y)-f(x)-\ev{G(x)}{y-x}}&\le M\om(\norm{y-x})\norm{y-x}
\end{align*}
for each $x,y\in E$;

\smallskip
\item[$\cW$]
if there exists a mapping $G\colon E\to X^*$ such that for each bounded $B\subset E$ the function $f$ is bounded on $B$, the mapping $G$ is bounded and uniformly continuous on $B$,
and for each $\ve>0$ there exists $\de>0$ such that
\[
\absb{big}{f(y)-f(x)-\ev{G(x)}{y-x}}\le\ve\norm{y-x}
\]
whenever $x,y\in B$ and $\norm{y-x}<\de$.
\end{itemize}
\end{definition}

Using only the continuity of $\om$ at $0$ it is easy to see that
\begin{equation}\label{e:Wom->WG}
\text{if $f$ satisfies $\cWG\om$ with some $G$, then it also satisfies $\cWd G$.}
\end{equation}
For some other connections between the above conditions see Lemmata~\ref{l:WG-concave}, \ref{l:WGlip-concave}, and Theorem~\ref{t:ext_C1ob}.
The motivation for the second inequality in $\cWG\om$ is the following fact which is an immediate consequence of the Taylor formula (see e.g.~\cite[Corollary~1.108]{HJ}).
\begin{lemma}\label{l:UTomega}
Let $U$ be an open convex subset of a normed linear space, $f\in C^1(U)$, and suppose that the derivative $Df$ is uniformly continuous with modulus $\om\in\mc M$.
Then
\[
\absb{big}{f(y)-f(x)-Df(x)[y-x]}\le\om(\norm{y-x})\norm{y-x}
\]
for each $x,y\in U$.
\end{lemma}
Indeed, Lemma~\ref{l:UTomega} clearly implies the following fact:
\begin{fact}\label{f:C1om->WGom}
Let $U$ be an open convex subset of a normed linear space, $\om\in\mc M$, and $f\in C^{1,\om}(U)$.
Then $f$ satisfies condition $\cWG\om$ (with $G=Df$).
\end{fact}
This fact will be generalised to quasiconvex sets in Lemma~\ref{l:quasic->WG}.

Sometimes we will work with renormed normed linear spaces,
so we need to know whether certain properties or notions (e.g. $C^{1,\om}$-smoothness) are invariant with respect to renormings:
\begin{lemma}\label{l:equiv}
Let $(X,\norm\cdot)$ be a normed linear space, $V\subset X$, $\om\in\mc M$, and let $\tnorm\cdot$ be an equivalent norm on $X$.
\begin{enumerate}[(i)]
\item
Assume that $\om$ is sub-additive or $V$ is convex.
Let $H\colon(V,\norm\cdot)\to(X^*,\norm\cdot^*)$ be uniformly continuous with modulus $\om$.
Then there is $K>0$ such that $H\colon(V,\tnorm\cdot)\to(X^*,\tnorm\cdot^*)$ is uniformly continuous with modulus $K\om$.
\item
Assume that $V$ is open, and $\om$ is sub-additive or $V$ is convex.
Then $C_{\norm\cdot}^{1,\om}(V)=C_{\tnorm\cdot}^{1,\om}(V)$.
\item
Assume that $V$ is open.
Then $C_{\norm\cdot}^{1,+}(V)=C_{\tnorm\cdot}^{1,+}(V)$.
\item
Assume that $\om$ is sub-additive and let $f\colon V\to\R$.
Then $f$ satisfies condition $\cWG\om$ on $V$ in the space $(X,\norm\cdot)$ if and only if $f$ satisfies $\cWG\om$ on $V$ in the space $(X,\tnorm\cdot)$ (with the same mapping $G$).
\end{enumerate}
\end{lemma}
\begin{proof}
Let $A\in\N$ be such that $\norm\cdot\le A\tnorm\cdot$.

(i)
Pick any $x,y\in V$.
If $\om$ is sub-additive, then by~\eqref{e:subad_Nt}
\[
\tnorm{H(x)-H(y)}^*\le A\norm{H(x)-H(y)}^*\le A\om(\norm{x-y})\le A\om(A\tnorm{x-y})\le A^2\om(\tnorm{x-y}).
\]
If $V$ is convex let us denote by $\om_H^{\norm\cdot}$ the minimal modulus of continuity of $H$ under the metrics given by $\norm\cdot$ and $\norm\cdot^*$.
Then $\om_H^{\norm\cdot}$ is sub-additive (Fact~\ref{f:conv_subad}), and therefore by the first part of the proof
\[
\tnorm{H(x)-H(y)}^*\le A^2\om_H^{\norm\cdot}(\tnorm{x-y})\le A^2\om(\tnorm{x-y}).
\]

(ii)
Note that by the symmetry it suffices to show just one inclusion.
Let $f\in C_{\norm\cdot}^{1,\om}(V)$ and let $L>0$ be such that $Df\colon(V,\norm\cdot)\to(X^*,\norm\cdot^*)$ is uniformly continuous with modulus $L\om$.
Then $Df\colon(V,\tnorm\cdot)\to(X^*,\tnorm\cdot^*)$ is uniformly continuous with modulus $KL\om$ by (i) and so $f\in C_{\tnorm\cdot}^{1,\om}(U)$.

(iii)
Note that by the symmetry it suffices to show just one inclusion.
Let $f\in C_{\norm\cdot}^{1,+}(V)$.
Choose any $\ve>0$.
There is $\de>0$ such that $\norm{Df(x)-Df(y)}^*<\frac\ve A$ whenever $x,y\in V$, $\norm{x-y}<\de$.
Now if $x,y\in V$ are such that $\tnorm{x-y}<\frac\de A$, then $\norm{x-y}\le A\tnorm{x-y}<\de$ and so $\tnorm{Df(x)-Df(y)}^*\le A\norm{Df(x)-Df(y)}^*<\ve$.

(iv)
By the symmetry it suffices to show just one implication.
So suppose that $f$ satisfies condition $\cWG\om$ on $V$ in the space $(X,\norm\cdot)$ with some $G\colon V\to X^*$ and $M>0$.
By (i) there is $K>0$ such that $\tnorm{G(x)-G(y)}^*\le KM\om(\tnorm{x-y})$ for each $x,y\in V$.
Further,
\[
\absb{big}{f(y)-f(x)-\ev{G(x)}{y-x}}\le M\om(\norm{y-x})\norm{y-x}\le AM\om(A\tnorm{y-x})\tnorm{y-x}\le A^2M\om(\tnorm{y-x})\tnorm{y-x}
\]
for each $x,y\in V$.
Thus $f$ satisfies $\cWG\om$ on $V$ in the space $(X,\tnorm\cdot)$ with $G$ and the constant $\wtilde M\eqdef\max\{A^2M, KM\}$.
\end{proof}

\subsection{\texorpdfstring{Semiconvex functions with modulus $\om$ and super-reflexive Banach spaces}{Semiconvex functions with modulus omega and super-reflexive Banach spaces}}

\begin{definition}
Let $X$ be a normed linear space, $U\subset X$ an open convex set, and $\om\in\mc M$.
We say that a function $f\colon U\to\R$ is semiconvex with modulus $\om$ (or $\om$-semiconvex for short) if
\[
f(\lambda x+(1-\lambda)y)\le\lambda f(x)+(1-\lambda)f(y)+\lambda(1-\lambda)\om(\norm{x-y})\norm{x-y}
\]
for every $x,y\in U$ and $\lambda\in[0,1]$.
We say that $f\colon U\to\R$ is $\om$-semiconcave if $-f$ is $\om$-semiconvex.
\end{definition}

\begin{remark}\label{r:paraconvex}
For the theory and applications of $\om$-semiconcave functions in $\R^n$ see~\cite{CS}.
For a comparison of $\om$-semiconvexity with the closely related Rollewicz's notion of strong paraconvexity
(which works with $\min\{\lambda,1-\lambda\}$ instead of $\lambda(1-\lambda)$), see e.g.~\cite[Remark 2.11]{DZ1} or~\cite[p.~218]{JTZ}.
\end{remark}

\begin{lemma}\label{l:inf-om_semiconvex}
Let $X$ be a normed linear space and $U\subset X$ open and convex.
Let $\{u_\al\}_{\al\in\Lambda}$ be a family of functions continuous on $U$ and semiconvex with the same modulus $\om\in\mc M$.
Then the function $u\eqdef\sup_{\al\in\Lambda}u_\al$ is also continuous and $\om$-semiconvex on $U$ provided that $u(x)<+\infty$ for each $x\in U$.
\end{lemma}
\begin{proof}
The $\om$-semiconvexity of $u$ easily follows from the definition of $\om$-semiconvexity (as in the proof of~\cite[Proposition~2.1.5]{CS}) and $u$ is clearly lower semi-continuous.
Since each lower semi-continuous and approximately convex (in the sense of~\cite{NLT}) function is continuous (even locally Lipschitz) by~\cite[Proposition 3.2]{NLT},
and $\om$-semiconvexity clearly implies approximate convexity, $u$ is continuous.
\end{proof}

For the following well-known fact see e.g.~\cite[(5), p.~838]{KZ} or \cite[Lemma~5.2]{DZ1}.
\begin{lemma}\label{l:C1om->semiconvex}
Let $X$ be a normed linear space, $U\subset X$ open convex set, and let $f\colon U\to\R$ be Fréchet differentiable such that $Df$ is uniformly continuous on $U$ with modulus $\om\in\mc M$.
Then $f$ is $\om$-semiconvex on $U$.
\end{lemma}

\begin{remark}
In the preceding lemma we can also clearly assert that $f$ is $\om$-semiconcave.
On the other hand if $f$ is continuous and both $\om$-semiconvex and $\om$-semiconcave, then for some $U$ (e.g. $U=X$) it follows that $f\in C^{1,\om}(U)$.
For more information see \cite{KZ} and \cite[Proposition~2.6]{AM}.
\end{remark}

The first basic ingredient for our paper is the following insertion (``in-between'') theorem which is a version of the well-known Ilmanen's lemma.

\begin{theorem}[{\cite[Corollary 3.2]{Kr1}}]\label{t:Ilmanen}
Let $X$ be a normed linear space, $\om\in\mc M$, and let $f_1\le f_2$ be continuous functions on $X$ such that $f_1$ is $\om$-semiconvex and $f_2$ is $\om$-semiconcave.
Then there exists $f\in C^{1,\om}(X)$ such that $f_1\le f\le f_2$.
\end{theorem}

\begin{remark}\hfill

(a)
Theorem~\ref{t:Ilmanen} is independently (implicitly and in a different language) proved in \cite[Proof of Theorem~3.2]{AM}.

(b)
Using \cite[Theorem~3.1]{Kr1} and \cite[Corollary~4.3]{KZ} we can assert that $Df$ is uniformly continuous with modulus $4\om$ in Theorem\ref{t:Ilmanen}.
\end{remark}

We now recall some facts we need about super-reflexive Banach spaces.
For the original definition of super-reflexive spaces and several of their characterisations see e.g.~\cite{DGZ} or~\cite{BL}.
For example a Banach space is super-reflexive if and only if it has an equivalent norm which is uniformly smooth.

Recall that if $X$ is a Banach space, then the \emph{modulus of smoothness} of the norm of $X$ (or of the space $X$) is defined as the function
\[
\rho(\tau)=\sup\bigg\{\frac{\norm{x+\tau y}+\norm{x-\tau y}}2-1\setsep x,y\in X, \norm x=\norm y=1\bigg\},\quad\tau\ge0.
\]
Notice that $\rho(\tau)\le\tau$ for any $\tau\ge0$.
We say that the modulus of smoothness is of power type $p>1$ if there exists $K>0$ such that $\rho(\tau)\le K\tau^p$ for each $\tau\ge0$; see e.g.~\cite{DGZ}.
Note that if $\rho$ is of power type $p>1$, then it is also of power type $q$ for any $1<q<p$ (recall that $\rho(\tau)\le\tau\le\tau^q$ for $\tau\ge1$).
For the following Pisier's result see e.g.~\cite{DGZ} or \cite[p.~412]{BL}.
\begin{theorem}\label{t:super-norm-modulus}
A Banach space is super-reflexive if and only if it has an equivalent norm with modulus of smoothness of power type $p$ for some $1<p\le2$.
\end{theorem}

\begin{remark}\label{r:Lp-modulus}
It is well-known (see~\cite[Corollary V.1.2]{DGZ}) that
if $\mu$ is an arbitrary measure, then $L_p(\mu)$ has modulus of smoothness of power type $p$ (resp.~$2$) if $1<p\le2$ (resp.~$p\ge2$).
In particular, every Hilbert space has modulus of smoothness of power type~$2$.
\end{remark}

We will need also the following well-known fact (see e.g.~\cite [Lemma 2.6]{DZ1} for an argument).
\begin{lemma}\label{l:norm_der_Holder}
If the norm $\norm\cdot$ of a Banach space $X$ has modulus of smoothness of power type $1+\al$ for some $0<\al\le1$,
then the Fréchet derivative $D\norm\cdot$ exists and is $\al$-Hölder on $S_X$.
\end{lemma}

The following lemma, which is the second basic ingredient of our paper, is a direct consequence of~\cite[Lemma 5.3]{DZ1} (together with Lemma~\ref{l:norm_der_Holder}).
\begin{lemma}\label{l:smooth_nu}
Let $\om\in\mc M$ be a concave modulus and let $X$ be a super-reflexive Banach space whose norm has modulus of smoothness of power type $1+\al$ for some $0<\al\le1$.
Denote $\vp(t)=\int_0^t\om(s)\d s$ for $t\ge0$ and $\nu(x)=\vp(\norm x)$ for $x\in X$.
Then the following assertions hold:
\begin{enumerate}[(i)]
\item
If $\om$ is bounded, then $\nu\in C^{1,\sigma}(X)$ for some bounded $\sigma\in\mc M$.
\item
If $\al=1$, then $\nu\in C^{1,\om}(X)$.
\item
If $\om(t)=ct^\beta$ for some $0<\beta\le\al$ and $c>0$, then $\nu\in C^{1,\om}(X)=C^{1,\beta}(X)$.
\end{enumerate}
\end{lemma}
Indeed, it is sufficient to use Lemma~\ref{l:norm_der_Holder} and then~\cite[Lemma 5.3]{DZ1} with $\beta\eqdef\al$, $\al\eqdef\beta$, $\vp\eqdef\om$, and $\psi\eqdef\vp$.
The items (i), (ii), (iii) follow from items (ii), (iv), (v) of~\cite[Lemma 5.3]{DZ1}, respectively.

\begin{remark}\label{r:relation_to_AM2}
Assertion (ii) and a slight generalisation of (iii) are essentially proved in~\cite{AM}.
Namely the proof of~\cite[Lemma 3.6]{AM} shows that $\nu\in C^{1,\om}(X)$ whenever the function $p(t)\eqdef t^\al/\om(t)$, $t>0$, is non-decreasing,
which clearly implies (iii) and also (ii), since for $\al=1$ the function $p$ is non-decreasing, because the concavity of $\om$ implies that $\om(t)/t=1/p(t)$ is non-increasing.
However, this latter fact is not mentioned explicitly in~\cite{AM}.
\end{remark}

\section{Basic lemmata and a proposition}\label{sec:basic}

The third basic ingredient of our paper is the following easy but non-trivial fact that is (without the explicit constant $6$) implicitly almost contained in~\cite{AM} without a proof or a reference
(cf. first two lines of p.~3 and the note just after Definition 1.4).

\begin{lemma}\label{l:WG->envelopes}
Let $X$ be a normed linear space, $E\subset X$, let $\om\in\mc M$ be concave, and suppose that $f\colon E\to\R$ satisfies condition $\cWG\om$ on $E$ with $G$ and $M>0$.
Let $\vp(t)=\int_0^t\om(s)\d s$ for $t\ge0$.
For each $z\in E$ set
\begin{align*}
h_z(x)&=f(z)+\ev{G(z)}{x-z}-6M\vp(\norm{x-z}),\\
H_z(x)&=f(z)+\ev{G(z)}{x-z}+6M\vp(\norm{x-z}),\quad x\in X,
\end{align*}
and define
\[
h(x)=\sup_{z\in E}h_z(x),\quad H(x)=\inf_{z\in E}H_z(x),\quad x\in X.
\]
Then $h(x)\le H(x)$ for each $x\in X$.
\end{lemma}
\begin{proof}
Let $x\in X$ and $z_1,z_2\in E$.
It is sufficient to show that $h_{z_1}(x)\le H_{z_2}(x)$, which clearly follows from
\[
\absb{big}{f(z_1)-f(z_2)+\ev{G(z_1)}{x-z_1}-\ev{G(z_2)}{x-z_2}}\le6M\vp(\norm{x-z_1})+6M\vp(\norm{x-z_2}).
\]
To prove this inequality we may assume that $\norm{x-z_1}\le\norm{x-z_2}$.
Using the concavity of $\om$ via Fact~\ref{f:concave-subad} and~Fact~\ref{f:om_vp} (twice) we obtain
\[\begin{split}
\om(\norm{z_1-z_2})\norm{z_1-z_2}&\le4\tfrac12\om\bigl(\tfrac12\norm{z_1-z_2}\bigr)\norm{z_1-z_2}\le8\vp\bigl(\tfrac12\norm{z_1-z_2}\bigr)\\
&\le8\vp\bigl(\tfrac12\norm{x-z_1}+\tfrac12\norm{x-z_2}\bigr)\le4\vp(\norm{x-z_1})+4\vp(\norm{x-z_2})
\end{split}\]
and using sub-additivity of $\om$ (Fact~\ref{f:concave-subad}) and Fact~\ref{f:om_vp} we obtain
\[\begin{split}
\om(\norm{z_1-z_2})\norm{x-z_1}&\le\bigl(\om(\norm{x-z_1})+\om(\norm{x-z_2})\bigr)\norm{x-z_1}\\
&\le\om(\norm{x-z_1})\norm{x-z_1}+\om(\norm{x-z_2})\norm{x-z_2}\le2\vp(\norm{x-z_1})+2\vp(\norm{x-z_2}).
\end{split}\]
Thus
\[\begin{split}
\absb{big}{f(z_1)-f(z_2)+\ev{G(z_1)}{x-z_1}-\ev{G(z_2)}{x-z_2}}&=\absb{big}{f(z_1)-f(z_2)-\ev{G(z_2)}{z_1-z_2}+\ev{G(z_1)-G(z_2)}{x-z_1}}\\
&\le M\om(\norm{z_1-z_2})\norm{z_1-z_2}+M\om(\norm{z_1-z_2})\norm{x-z_1}\\
&\le6M\vp(\norm{x-z_1})+6M\vp(\norm{x-z_2}).
\end{split}\]
\end{proof}

All our extension results are based on the proposition below (which is a consequence of Lemma~\ref{l:WG->envelopes} and Theorem~\ref{t:Ilmanen}) and Lemma~\ref{l:smooth_nu}.
This proposition with $\sigma=\om$ is essentially mentioned in~\cite{AM} (see the first sentence after \cite[Theorem~1.10]{AM}) and is implicitly used in~\cite{AM}.
The case $\sigma\neq\om$ together with Lemma~\ref{l:smooth_nu}(i) will be used substantially in Section~\ref{sec:C1+}.

\begin{proposition}\label{p:WGom->extension}
Let $\om\in\mc M$ be concave, $\vp(t)=\int_0^t\om(s)\d s$ for $t\ge0$, and $\sigma\in\mc M$.
Let $X$ be a normed linear space such that the function $\nu(x)=\vp(\norm x)$, $x\in X$, is $C^{1,\sigma}$-smooth.
Let $f$ be a real function on $E\subset X$ which satisfies condition $\cWG\om$ on $E$ with $G\colon E\to X^*$.
Then $f$ can be extended to a function $F\in C^{1,\sigma}(X)$ such that $DF(x)=G(x)$ for each $x\in E$.
\end{proposition}
\begin{proof}
Let $M>0$ be as in condition $\cWG\om$.
For each $z\in E$ set
\begin{align*}
h_z(x)&=f(z)+\ev{G(z)}{x-z}-6M\vp(\norm{x-z}),\\
H_z(x)&=f(z)+\ev{G(z)}{x-z}+6M\vp(\norm{x-z}),\quad x\in X,
\end{align*}
and define
\[
h(x)=\sup_{z\in E}h_z(x),\quad H(x)=\inf_{z\in E}H_z(x),\quad x\in X.
\]
Then $h\le H$ by Lemma~\ref{l:WG->envelopes}.
Let $x\in E$.
Since $h_x(x)=H_x(x)=f(x)$, it follows that $H(x)\le f(x)\le h(x)$ and consequently $h(x)=f(x)=H(x)$.

By the assumptions there exists $K>0$ such that $D\nu$ is uniformly continuous with modulus $K\sigma$,
which easily implies that for each $z\in E$ the derivative $Dh_z$ is uniformly continuous with modulus $6MK\sigma$.
So each function $h_z$, $z\in E$, is continuous and semiconvex with modulus $6MK\sigma$ by Lemma~\ref{l:C1om->semiconvex}
and so also (clearly finite) function $h$ is continuous and semiconvex with modulus $6MK\sigma$ by Lemma~\ref{l:inf-om_semiconvex}.
Quite analogously (working with functions $-H_z$ and $-H=\sup_{z\in E}-H_z$) we obtain that $H$ is semiconcave with modulus $6MK\sigma$ and continuous.
So Theorem~\ref{t:Ilmanen} implies that there exists $F\in C^{1,6MK\sigma}(X)=C^{1,\sigma}(X)$ such that $h\le F\le H$.
Since $h(x)=H(x)=f(x)$ for $x\in E$, it follows that $F$ is an extension of $f$.

Further, fix $z\in E$.
Then $h_z(x)\le h(x)\le F(x)\le H(x)\le H_z(x)$ for each $x\in X$.
Since $D\nu(0)$ exists by the assumption and $\nu$ has a minimum at~$0$, it is clear that $D\nu(0)=0$.
Therefore $Dh_z(z)=G(z)=DH_z(z)$.
Since also $h_z\le F\le H_z$ and $h_z(z)=F(z)=H_z(z)$, it is easy to see that $DF(z)=G(z)$.
\end{proof}

\begin{remark}
The smoothness of $\nu$ and the introduction of functions $h_z$ and $h$ (with $M=1/6$) was essential already in~\cite{DZ1}
in connection with extensions of $\om$-semiconvex functions; see~\cite[Lemma 5.4 and Proposition 5.12]{DZ1}.
\end{remark}

The following fact is a basic tool for our extension results from open quasiconvex sets.

\begin{lemma}\label{l:quasic->WG}
Let $U$ be an open $c$-quasiconvex subset of a normed linear space $X$, $\om\in\mc M$, and $f\in C^{1,\om}(U)$.
Then $f$ satisfies condition $\cWG\om$ with $G=Df$ (and $M=4Kc^3$, where $K>0$ is such that $\om_{Df}\le K\om$).
\end{lemma}
\begin{proof}
Consider an arbitrary $K>0$ such that $\om_{Df}\le K\om$.
By Lemma~\ref{l:quasic-concave_modulus} there is a concave $\wtilde\om\in\mc M$ such that $\frac12\wtilde\om\le\om_{Df}\le c\wtilde\om$.
Note that $\wtilde\om\le2\om_{Df}\le2K\om$.
Consider arbitrary different points $x,y\in U$.
Then we can choose a continuous $\ga\colon[0,1]\to U$ such that $\ga(0)=x$, $\ga(1)=y$, and $0<\norm{x-y}\le L\eqdef\len\ga\le c\norm{x-y}$.
Set $\ve=\dist(\ga([0,1]),X\setminus U)$.
Then $\ve>0$ by the compactness of $\ga([0,1])$.
By the uniform continuity of $\ga$ we choose $\de>0$ such that $\norm{\ga(s)-\ga(t)}<\ve$ whenever $0\le s\le t\le1$ and $t-s<\de$.
Further, choose points $0=t_0<t_1<\dotsb<t_{n-1}<t_n=1$ such that $t_i-t_{i-1}<\de$ for $i=1,\dotsc,n$ and denote $x_i=\ga(t_i)$.
Then clearly $x_0=x$, $x_n=y$, and
\begin{equation}\label{e:blizko}
\norm{x_j-x_0}\le\sum_{i=1}^n\norm{x_i-x_{i-1}}\le L,\quad j=1,\dotsc,n.
\end{equation}
The choice of $\ve$, $\de$, and $t_0,\dotsc,t_n$ implies that $x_i\in U(x_{i-1},\ve)\subset U$ for $i=1,\dotsc,n$,
and so Lemma~\ref{l:UTomega} together with $\om_{Df}\le c\wtilde\om$ and~\eqref{e:blizko} imply that
\[
\absb{big}{f(x_i)-f(x_{i-1})-Df(x_{i-1})[x_i-x_{i-1}]}\le c\wtilde\om(\norm{x_i- x_{i-1}})\norm{x_i-x_{i-1}}\le c\wtilde\om(L)\norm{x_i-x_{i-1}}
\]
for $i=1,\dotsc,n$.
Further, $\om_{Df}\le c\wtilde\om$ with~\eqref{e:blizko} gives $\norm{Df(x_j)-Df(x_0)}\le c\wtilde\om(\norm{x_j-x_0})\le c\wtilde\om(L)$ for $j=1,\dotsc,n$.
Using these inequalities, \eqref{e:blizko}, the concavity of $\wtilde\om$, and $c\ge1$ we obtain
\[\begin{split}
\absb{big}{f(y)-f(x)-Df(x)[y-x]}&=\absa{\sum_{i=1}^n\bigl(f(x_i)-f(x_{i-1})\bigr)-\sum_{i=1}^nDf(x_0)[x_i-x_{i-1}]}\\
&\le\absa{\sum_{i=1}^n\bigl(f(x_i)-f(x_{i-1})-Df(x_{i-1})[x_i-x_{i-1}]\bigr)}+\absa{\sum_{i=1}^n\bigl(Df(x_{i-1})-Df(x_0)\bigr)[x_i-x_{i-1}]}\\
&\le\sum_{i=1}^nc\wtilde\om(L)\norm{x_i-x_{i-1}}+\sum_{i=1}^nc\wtilde\om(L)\norm{x_i-x_{i-1}}\le2c\wtilde\om(L)L\\
&\le2c^2\wtilde\om(c\norm{y-x})\norm{y-x}\le2c^3\wtilde\om(\norm{y-x})\norm{y-x}\le4Kc^3\om(\norm{y-x})\norm{y-x},
\end{split}\]
and so $f$ satisfies $\cWG\om$ with $G=Df$ and $M=4Kc^3$, since $\om_{Df}\le K\om\le M\om$.
\end{proof}

\section{\texorpdfstring{The $C^{1,\om}$ case}{The C1omega case}}\label{sec:C1om}

The main result of this section is Theorem~\ref{t:ext_C1om}.
For the convenience of the reader we will repeat the statement.
\begin{th_ex_om}
Let $\om\in\mc M$ be a concave modulus and let $X$ be a super-reflexive Banach space that has an equivalent norm with modulus of smoothness of power type $2$.
Let $E\subset X$ and let $f$ be a real function on $E$.
Then $f$ can be extended to a function $F\in C^{1,\om}(X)$ if and only if $f$ satisfies condition $\cWG\om$ on $E$.
Moreover, if $\cWG\om$ is satisfied, then $F$ can be found such that $DF(x)=G(x)$ for each $x\in E$.
\end{th_ex_om}
The assumption on $X$ cannot be relaxed, see Remark~\ref{r:part_unity}(a).
\begin{proof}
$\Rightarrow$
This implication holds in fact in an arbitrary normed linear space~$X$.
Indeed, Fact~\ref{f:C1om->WGom} implies that $F$ satisfies $\cWG\om$ on $X$ and so $f$ satisfies $\cWG\om$ on $E$.

$\Leftarrow$
Choose an equivalent norm $\tnorm\cdot$ on $X$ with modulus of smoothness of power type~$2$.
Then $f$ satisfies condition $\cWG\om$ with respect to the norm $\tnorm\cdot$ with some $G$ (Lemma~\ref{l:equiv}(iv)).
The function $\nu=\vp\comp\tnorm\cdot$, where $\vp(t)=\int_0^t\om(s)\d s$, is $C^{1,\om}$-smooth on $(X,\tnorm\cdot)$ by Lemma~\ref{l:smooth_nu}(ii).
So we can apply Proposition~\ref{p:WGom->extension} with $\sigma=\om$ and extend $f$ to a function $F$ on $X$ which is $C^{1,\om}$-smooth on $(X,\tnorm\cdot)$,
and hence also on $(X,\norm\cdot)$ by Lemma~\ref{l:equiv}(ii).
Moreover, Proposition~\ref{p:WGom->extension} gives that $DF(x)=G(x)$ for each $x\in E$.
\end{proof}

As a corollary we obtain the following true extension theorem.
\begin{corollary}\label{c:extC1om-quasic}
Let $X$ be a super-reflexive Banach space that has an equivalent norm with modulus of smoothness of power type $2$.
Let $U\subset X$ be an open quasiconvex set, $\om\in\mc M$, and $f\in C^{1,\om}(U)$.
Then $f$ can be extended to a function $F\in C^{1,\om}(X)$.
\end{corollary}
The assumption on $X$ cannot be relaxed, see Remark~\ref{r:quasic-bump}.
\begin{proof}
There is $K>0$ such that $\om_{Df}\le K\om$.
By Lemma~\ref{l:quasic-concave_modulus} there are a concave $\wtilde\om\in\mc M$ and $c\ge1$ such that $\frac12\wtilde\om\le\om_{Df}\le c\wtilde\om$.
Therefore $f\in C^{1,\wtilde\om}(U)$ and Lemma~\ref{l:quasic->WG} then implies that $f$ satisfies condition $\cWG{\wtilde\om}$ on $U$.
Thus $f$ can be extended to a function $F\in C^{1,\wtilde\om}(X)$ by Theorem~\ref{t:ext_C1om}.
Since $\wtilde\om\le2\om_{Df}\le2K\om$, it follows that $F\in C^{1,\om}(X)$.
\end{proof}

Corollary~\ref{c:extC1om-quasic} for $X=\R^n$ and $\om(t)=t$, $t\ge0$, ``almost follows'' from~\cite[Theorem 2.64]{BB}; for the details see the following remark.
\begin{remark}\label{r:brud}
(a)
\cite[Theorem 2.64]{BB} gives that if $\om\in\mc M$ is concave, $U\subset\R^n$ is a $(C,\om)$-convex domain, $k\ge1$, and $f\in\dot C^{k,\om}(U)$, then $f$ has an extension $F\in\dot C^{k,\om}(\R^n)$.
Note that if $f\in C^{1,\om}(U)$ and both $f$ and $Df$ are bounded, then clearly $f\in\dot C^{1,\om}(U)$ (and the opposite implication holds for $U=\R^n$).
Further, if $\om(t)=t$ for $t\ge0$, then $(C,\om)$-convexity coincides with $C$-quasiconvexity.
Using these facts we obtain that the assertion of Corollary~\ref{c:extC1om-quasic} follows from \cite[Theorem 2.64]{BB} if $X=\R^n$, $\om(t)=t$ for $t\ge0$, and both $f$ and $Df$ are bounded on $U$.

(b)
However, this connection between Corollary~\ref{c:extC1om-quasic} and \cite[Theorem 2.64]{BB} fails for some more general moduli.
For example, if $0<\al<1$ is fixed and $\om_\al(t)=t^\al$ for $t\ge0$,
then the set of all open quasiconvex sets in $\R^n$ ($n>1$) is strictly larger than the set of all $(C,\om_\al)$-convex domains with any $C>0$.
Indeed, suppose that $U\subset\R^n$ is a $(C,\om_\al)$-convex domain, i.e.
for each $x,y\in U$ there exists a ``polygonal line $\ga\colon[0,1]\to U$ with the segments $[\ga(t_i),\ga(t_{i+1})]$, where $0=t_0<t_1<\dotsb<t_k=1$'' such that
\[
\ell_{\om_\al}(\ga)\eqdef\sum_{i=0}^{k-1}\norm{\ga(t_{i+1})-\ga(t_i)}^\al\le C\norm{x-y}^\al.
\]
Then
\[
\len\ga=\sum_{i=0}^{k-1}\norm{\ga(t_{i+1})-\ga(t_i)}\le\left(\sum_{i=0}^{k-1}\norm{\ga(t_{i+1})-\ga(t_i)}^\al\right)^{\frac1\al}\le C^{\frac1\al}\norm{x-y}.
\]
(For the first inequality recall that $\om_\al$ is sub-additive.)
Therefore $U$ is $C^{\frac1\al}$-quasiconvex.
Note that this fact that we just proved contradicts the incorrect claim from \cite{BB} that the set $G_{\frac1\al}$ from \cite[Figure 2.1]{BB} (cf. Remark~\ref{r:ext-lip-quasic}), which is clearly not quasiconvex, is $(C,\om_\al)$-convex.
For another argument that this claim is incorrect see Remark~\ref{r:ext-lip-quasic}.

Further, in the case $\al=\frac12$ (other cases are similar) we can define an open quasiconvex set $U\subset\R^2$ which is not a $(C,\om_\al)$-convex domain for any $C$ by
\[
U=\left\{(x,y)\in\R^2\setsep 0<x<1, x^2\sin\frac1x<y<x^4+x^2\sin\frac1x\right\}.
\]
It is not too difficult to prove that $U$ has the aforementioned properties.
\end{remark}

\begin{remark}\label{r:ext-quasic}
It is well-known that Corollary \ref{c:extC1om-quasic} (resp. Corollary~\ref{c:extC1_Hold-quasic} below)
does not hold for general simply connected domains~$U$ (i.e. we cannot relax the assumption of quasiconvexity of $U$ to being simply connected) even when $X=\R^2$ and $\om\in\mc M$ is arbitrary concave non-zero (resp. $0<\al\le1$).
Since we were not able to find a reference to a corresponding construction, we supply the following easy one:
Let $U=(0,3)\times(0,3)\setminus\bigl([1,2]\times[1,2]\cup(0,1)\times\{2\}\bigr)\subset\R^2$ and define $f\colon U\to\R$ by
\[
f(x,y)=\begin{cases}
0\quad&\text{for $(x,y)\in(0,3)\times(0,1)\cup(0,1)\times[1,2)$,}\\
1\quad&\text{for $(x,y)\in(0,3)\times(2,3)$,}\\
\sin^2\bigl(\tfrac\pi2(y-1)\bigr)\quad&\text{for $(x,y)\in(2,3)\times[1,2]$.}
\end{cases}
\]
Then $f\in C^{1,1}(U)\subset C^{1,\om}(U)$, but it cannot be extended even to a continuous function on $\cl[1]U$.
\end{remark}

Now we prove a proposition which works with arbitrary (possibly \emph{non-concave}) moduli $\om\in\mc M$ and may be regarded as a generalisation of Theorem~\ref{t:ext_C1om}.

\begin{proposition}\label{p:ext_C1om}
Let $X$ be a super-reflexive Banach space that has an equivalent norm with modulus of smoothness of power type~$2$.
Let $E\subset X$, let $f$ be a real function on $E$, and let $\om\in\mc M$.
Then $f$ can be extended to a function $F\in C^{1,\om}(X)$ if and only if there exists a concave $\psi\in\mc M$ such that $\psi\le\om$ and $f$ satisfies $\cWG\psi$ on~$E$.
Moreover, if $f$ satisfies $\cWG\psi$ with some $G$, then $F$ can be found such that $DF(x)=G(x)$ for each $x\in E$.
\end{proposition}
\begin{proof}
$\Rightarrow$
Let $K>0$ be such that $\om_{DF}\le K\om$.
The modulus $\om_{DF}$ is sub-additive (Fact~\ref{f:conv_subad}) and so there is a concave $\wtilde\om\in\mc M$ such that $\om_{DF}\le\wtilde\om\le2\om_{DF}\le2K\om$ (Lemma~\ref{l:Steckin}(iii)).
Thus we can define concave modulus $\psi=\frac1{2K}\wtilde\om$ for which $\psi\le\om$ and $F\in C^{1,\psi}(X)$.
Fact~\ref{f:C1om->WGom} implies that $F$ satisfies $\cWG\psi$ on $X$ and so $f$ satisfies $\cWG\psi$ on $E$.

$\Leftarrow$
By Theorem~\ref{t:ext_C1om} the function $f$ can be extended to $F\in C^{1,\psi}(X)\subset C^{1,\om}(X)$ and the ``moreover'' part holds.
\end{proof}

In connection with the previous proposition, in the following remark we consider the validity of Theorem~\ref{t:ext_C1om} for non-concave moduli.
The main observation is that the assumption of concavity cannot be completely dropped.
\begin{remark}\label{r:moduli_nonconc}
(a)
Assume that $\om\in\mc M$ is not concave and Theorem~\ref{t:ext_C1om} holds for this~$\om$.
Then there is a non-zero concave $\psi\in\mc M$ such that $\psi\le\om$ (which is clearly equivalent to the fact that $\liminf_{t\to0+}\frac{\om(t)}t>0$).
Indeed, consider $e\in X$, $e\neq0$.
Set $f(0)=0$, $f(e)=0$, $f(2e)=1$.
Then $f$ clearly satisfies $\cWG\om$ (e.g. with $G=0$) and so can be extended to $F_1\in C^{1,\om}(X)$.
By Proposition~\ref{p:ext_C1om} there exists a concave $\psi\in\mc M$ such that $\psi\le\om$ and $f$ satisfies $\cWG\psi$.
Thus by Theorem~\ref{t:ext_C1om} the function $f$ can be extended to $F_2\in C^{1,\psi}(X)$.
It follows that $\psi$ is non-zero.
Otherwise $DF_2$ would be constant and in turn $F_2$ would be affine, which is a contradiction with the definition of~$f$.

So we have shown that Theorem~\ref{t:ext_C1om} does not hold for arbitrary $\om\in\mc M$.

(b)
Of course there are some non-concave moduli $\om$ for which Theorem~\ref{t:ext_C1om} holds.
In fact it is easy to see that it is sufficient to assume that there exist a concave $\psi\in\mc M$ and $K>0$ such that $\psi\le\om\le K\psi$, which holds e.g. if $\om$ is sub-additive (Lemma~\ref{l:Steckin}(iii)).
However, we do not know whether this condition is necessary.
\end{remark}

Other interesting consequences of Proposition~\ref{p:WGom->extension} are the following extension results for functions with Hölder derivatives.

\begin{theorem}\label{t:extC1_Holder}
Let $0<\al\le1$ and let $X$ be a super-reflexive Banach space that has an equivalent norm with modulus of smoothness of power type $1+\al$.
Let $E\subset X$ and let $f$ be a real function on $E$.
Then $f$ can be extended to a function $F\in C^{1,\al}(X)$ if and only if condition $\cWG\om$ holds for $\om(t)=t^{\al}$.

Moreover, if $\cWG\om$ is satisfied, then $F$ can be found such that $DF(x)=G(x)$ for $x\in E$.
\end{theorem}
The assumption on $X$ cannot be relaxed, see Remark~\ref{r:part_unity}(a).
\begin{proof}
$\Rightarrow$ follows from Fact~\ref{f:C1om->WGom} in an arbitrary normed linear space as in the proof of Theorem~\ref{t:ext_C1om}.

$\Leftarrow$
Choose an equivalent norm $\tnorm\cdot$ on $X$ with modulus of smoothness of power type~$1+\al$.
Then $f$ satisfies condition $\cWG\om$ with respect to the norm $\tnorm\cdot$ with some $G$ (Lemma~\ref{l:equiv}(iv)).
The function $\nu=\vp\comp\tnorm\cdot$, where $\vp(t)=\int_0^t\om(s)\d s$, is $C^{1,\om}$-smooth on $(X,\tnorm\cdot)$ by Lemma~\ref{l:smooth_nu}(iii).
So we can apply Proposition~\ref{p:WGom->extension} with $\sigma=\om$ and extend $f$ to a function $F$ on $X$ which is $C^{1,\om}$-smooth on $(X,\tnorm\cdot)$,
and hence also on $(X,\norm\cdot)$ by Lemma~\ref{l:equiv}(ii).
Moreover, Proposition~\ref{p:WGom->extension} gives that $DF(x)=G(x)$ for each $x\in E$.
\end{proof}

By Remark~\ref{r:Lp-modulus}, Theorem~\ref{t:extC1_Holder} (and its Corollary~\ref{c:extC1_Hold-quasic}) can be applied e.g. if $X$ is isomorphic to some $L_p(\mu)$ with $p\ge1+\al$.

\begin{corollary}\label{c:extC1_Hold-quasic}
Let $0<\al\le1$ and let $X$ be a super-reflexive Banach space that has an equivalent norm with modulus of smoothness of power type $1+\al$.
Let $U\subset X$ be an open quasiconvex set and $f\in C^{1,\al}(U)$.
Then $f$ can be extended to a function $F\in C^{1,\al}(X)$.
\end{corollary}
The assumption on $X$ cannot be relaxed, see Remark~\ref{r:quasic-bump}.
\begin{proof}
Set $\om(t)=t^\al$ for $t\ge0$.
Lemma~\ref{l:quasic->WG} implies that $f$ satisfies condition $\cWG\om$ on $U$ and so $f$ can be extended to a function $F\in C^{1,\al}(X)$ by Theorem~\ref{t:extC1_Holder}.
\end{proof}

\begin{remark}\label{r:part_unity}
Let $X$ be a normed linear space and let $\om\in\mc M$ be such that $\om(t)>0$ for $t>0$.
Assume that we can extend every function from an arbitrary closed subset $E$ of $X$ satisfying condition $\cWG\om$ on $E$ to a function in $C^{1,\om}(X)$.
Then we obtain the following smooth separation property:
For every closed $A\subset X$ and $\de>0$ there is a function $\vp\in C^{1,\om}(X)$ such that $\vp=1$ on $A$ and $\vp(x)=0$ whenever $\dist(x,A)\ge\de$.
Indeed, set $f=1$, $G=0$ on $A$ and $f=0$, $G=0$ on $B=\{x\in X\setsep\dist(x,A)\ge\de\}$.
Then $f$ satisfies condition $\cWG\om$ with $G$ and $M=\frac1{\de\om(\de)}$ on the closed set $A\cup B$ and so it can be extended to a desired function $\vp\in C^{1,\om}(X)$.
In particular, we can produce a $C^{1,\om}$-smooth bump function on~$X$ and consequently $X$ is super-reflexive by \cite[Theorem~V.3.2]{DGZ}.
(A similar observation is already known, cf. \cite[p.~9]{AM}.)

In the Hölder case, i.e. $\om(t)=t^\al$ for some $0<\al\le1$, we obtain two new interesting applications.

(a)
The assumption on $X$ in Theorems~\ref{t:ext_C1om} and~\ref{t:extC1_Holder} (i.e. the existence of an equivalent norm with modulus of smoothness of power type $2$, resp. $1+\al$) cannot be relaxed
(for Theorem~\ref{t:ext_C1om} consider $\om(t)=t$).
This follows again from \cite[Theorem~V.3.2]{DGZ}.

(b)
The above separation property gives the existence of a locally finite $C^{1,\al}$-smooth partition of unity on an arbitrary open $U\subset X$ subordinated to an arbitrary open covering of~$U$.
Indeed, we use \cite[Lemma~7.49, (ii)$\Rightarrow$(vi)]{HJ} with $S=\{\vp\in C^{1,\al}(U)\setsep\text{$\vp$ is bounded and Lipschitz}\}$.
Then $S$ is a partition ring, which can be seen by repeating the beginning of the proof of \cite[Theorem~7.56]{HJ} with obvious changes.
To verify condition \cite[Lemma~7.49(ii)]{HJ} we set $A=\cl[1]V$, $\de=\dist(V,U\setminus W)$ and use the separation property above.
This gives us a corresponding $\psi\in C^{1,\al}(X)$ and since it has a bounded support, both $\psi$ and its derivative are bounded and so $\psi$ is Lipschitz.
Then $\vp=\psi\restr U\in S$.

Therefore as a corollary of Theorem~\ref{t:extC1_Holder} we obtain the existence of locally finite $C^{1,\al}$-smooth partitions of unity on (open subsets of) super-reflexive spaces that have an equivalent norm with modulus of smoothness of power type $1+\al$.
The proof we just described is much easier than the original proof in~\cite{JTZ:UFPartUnity} which uses some deep and heavy machinery of non-separable Banach space theory.
So we believe that the above observation is interesting in itself.
\end{remark}

\section{\texorpdfstring{The $C^{1,+}$ case}{The C1+ case}}\label{sec:C1+}

First we observe two more connections (deeper than \eqref{e:Wom->WG}) between the Whitney-Glaeser type conditions from Definition~\ref{d:condW}.

\begin{lemma}\label{l:WG-concave}
Let $X$ be a normed linear space, let $f$ be a function on $E\subset X$, and let $G\colon E\to X^*$.
Then the following statements are equivalent:
\begin{enumerate}[(i)]
\item $f$ satisfies condition $\cWd G$ and for each $d>0$ there exists $K>0$ such that
\begin{align*}
\norm{G(y)-G(x)}&\le K\norm{y-x},\\
\absb{big}{f(y)-f(x)-\ev{G(x)}{y-x}}&\le K\norm{y-x}^2
\end{align*}
whenever $x,y\in E$, $\norm{y-x}\ge d$.
\item
There exists a concave modulus $\om\in\mc M$ such that $f$ satisfies condition $\cWG\om$ with $G$ on~$E$.
\end{enumerate}
\end{lemma}
\begin{proof}
(ii)$\Rightarrow$(i)
Let $M>0$ be the constant from condition $\cWG\om$.
The function $f$ clearly satisfies condition $\cWd G$ (see \eqref{e:Wom->WG}).
Now let $d>0$ be given.
Then $\om(t)=\om\bigl(\frac tdd\bigr)\le\frac td\om(d)$ for $t\ge d$ by Fact~\ref{f:concave-subad}.
It follows that
\begin{align*}
\norm{G(y)-G(x)}&\le M\om(\norm{y-x})\le M\frac{\om(d)}d\norm{y-x},\\
\absb{big}{f(y)-f(x)-\ev{G(x)}{y-x}}&\le M\om(\norm{y-x})\norm{y-x}\le M\frac{\om(d)}d\norm{y-x}^2
\end{align*}
whenever $x,y\in E$, $\norm{y-x}\ge d$.

(i)$\Rightarrow$(ii)
For each $\de\ge0$ set
\[
\al(\de)=\!\!\sup_{\substack{x,y\in E\\0<\norm{y-x}\le\de}}\!\!\max\left\{\norm{G(y)-G(x)},\frac{\absb{big}{f(y)-f(x)-\ev{G(x)}{y-x}}}{\norm{y-x}}\right\}
\]
if $\{(x,y)\in E\times E\setsep0<\norm{y-x}\le\de\}\neq\emptyset$ and $\al(\de)=0$ otherwise.
Obviously, $\al\colon[0,+\infty)\to[0,+\infty]$ is non-decreasing and $\al(0)=0$.
Using $\cWd G$ we easily see that $\lim_{\de\to 0_+}\al(\de)=0$.
In particular there exists $\de_0>0$ such that $\al(\de)<1$ for each $0<\de\le\de_0$.
Let $K>0$ be the constant from (i) corresponding to $d=\de_0$.
Then
\[
\norm{G(y)-G(x)}\le K\norm{y-x}\le K\de\quad\text{and}\quad\frac{\absb{big}{f(y)-f(x)-\ev{G(x)}{y-x}}}{\norm{y-x}}\le K\norm{y-x}\le K\de
\]
whenever $x,y\in E$ and $\de_0\le\norm{y-x}\le\de$, and so $\al(\de)\le\max\{1,K\de\}$ for each $\de\ge0$.
Hence $\al\in\mc M$ and by Lemma~\ref{l:Steckin}(i) there exists a concave modulus $\om\in\mc M$ such that $\al\le\om$.
Now if arbitrary $u,v\in E$, $u\neq v$, are given, then
$\max\bigl\{\norm{G(v)-G(u)},\frac{\abs{f(v)-f(u)-\ev{G(u)}{v-u}}}{\norm{v-u}}\bigr\}\le\al(\norm{v-u})\le\om(\norm{v-u})$
and so $f$ satisfies $\cWG\om$ on $E$ with $G$ and $M=1$.
\end{proof}

\begin{lemma}\label{l:WGlip-concave}
Let $X$ be a normed linear space, let $f$ be a function on $E\subset X$, and let $G\colon E\to X^*$.
Then the following statements are equivalent:
\begin{enumerate}[(i)]
\item $f$ is Lipschitz, $G$ is bounded, and $f$ satisfies condition $\cWd G$.
\item There exists a bounded concave modulus $\om\in\mc M$ such that $f$ satisfies condition $\cWG\om$ with $G$ on~$E$.
\end{enumerate}
\end{lemma}
\begin{proof}
(ii)$\Rightarrow$(i)
The function $f$ clearly satisfies condition $\cWd G$ (see \eqref{e:Wom->WG}).
Further, suppose that $f$ satisfies condition $\cWG\om$ on~$E$ with $G$ and some $M>0$.
Let $K\ge0$ be such that $\om(t)\le K$ for $t\ge0$.
Fix $x_0\in E$ and set $L=\norm{G(x_0)}+MK$.
Then $\norm{G(x)}\le\norm{G(x_0)}+\norm{G(x)-G(x_0)}\le\norm{G(x_0)}+M\om(\norm{x-x_0})\le L$ for any $x\in E$.
Finally,
\[
\abs{f(y)-f(x)}\le\absb{big}{f(y)-f(x)-\ev{G(x)}{y-x}}+\abs{\ev{G(x)}{y-x}}\le M\om(\norm{y-x})\norm{y-x}+L\norm{y-x}\le(KM+L)\norm{y-x}
\]
for any $x,y\in E$.

(i)$\Rightarrow$(ii)
Let $L>0$ be such that $f$ is $L$-Lipschitz and $\norm{G(x)}\le L$ for each $x\in E$.
Then
\begin{equation}\label{e:WG-lip}
\norm{G(y)-G(x)}\le2L\quad\text{and}\quad\absb{big}{f(y)-f(x)-\ev{G(x)}{y-x}}\le2L\norm{y-x}
\end{equation}
whenever $x,y\in E$.
So, if $d>0$ and $x,y\in E$ are such that $\norm{y-x}\ge d$, then $\norm{G(y)-G(x)}\le\frac{2L}d\norm{y-x}$ and $\absb{big}{f(y)-f(x)-\ev{G(x)}{y-x}}\le\frac{2L}d\norm{y-x}^2$.
Therefore by Lemma~\ref{l:WG-concave} there is a concave $\wtilde\om\in\mc M$ such that $f$ satisfies condition $\cWG{\wtilde\om}$ on~$E$ with $G$ and a constant $M>0$.
Using this fact and~\eqref{e:WG-lip} we obtain that $f$ satisfies condition $\cWG\om$ on~$E$ with $G$ and the constant $M$, where $\om=\min\bigl\{\wtilde\om,\frac{2L}M\bigr\}$.
\end{proof}

\begin{theorem}\label{t:ext_C1+_gen}
Let $X$ be a normed linear space, $E\subset X$, and let $f$ be a real function on $E$.
Consider the following statements:
\begin{enumerate}[(i)]
\item $f$ can be extended to a function $F\in C^{1,+}(X)$.
\item There exists $G$ such that $f$ satisfies condition $\cWd G$ and for each $d>0$ there exists $K>0$ such that
\begin{align*}
\norm{G(y)-G(x)}&\le K\norm{y-x},\\
\absb{big}{f(y)-f(x)-\ev{G(x)}{y-x}}&\le K\norm{y-x}^2
\end{align*}
whenever $x,y\in E$, $\norm{y-x}\ge d$.
\end{enumerate}
Then (i)$\Rightarrow$(ii).
If $X$ is a super-reflexive Banach space that has an equivalent norm with modulus of smoothness of power type~$2$, then both statements are equivalent
and moreover if (ii) holds, then $F$ can be found such that $DF(x)=G(x)$ for each $x\in E$.
\end{theorem}
\begin{proof}
(i)$\Rightarrow$(ii)
By Fact~\ref{f:conv_subad}, $\om_{DF}\in\mc M$ is sub-additive and so by Lemma~\ref{l:Steckin}(iii) there is a concave $\om\in\mc M$ such that $\om_{DF}\le\om$.
Fact~\ref{f:C1om->WGom} used on $F$ implies that $f$ satisfies condition $\cWG\om$ on $E$ with $G=DF\restr E$ and we may apply Lemma~\ref{l:WG-concave}.

(ii)$\Rightarrow$(i)
By Lemma~\ref{l:WG-concave} there exists a concave modulus $\om\in\mc M$ such that $f$ satisfies condition $\cWG\om$ with $G$ on~$E$.
By Theorem~\ref{t:ext_C1om} the function $f$ can be extended to a function $F\in C^{1,\om}(X)\subset C^{1,+}(X)$ such that $DF(x)=G(x)$ for each $x\in E$.
\end{proof}

\begin{remark}\label{r:ext-gen}
We do not know whether the assumption on~$X$ for (ii)$\Rightarrow$(i) in the preceding theorem (i.e. the existence of an equivalent norm with modulus of smoothness of power type~$2$) can be relaxed.
However, the super-reflexivity of $X$ is necessary:
An extension of the function defined by $f=0$ on $X\setminus B_X$ and $f(0)=1$ (which satisfies (ii) with $G=0$) produces a $C^{1,+}$-smooth bump,
so we can use \cite[Theorem~V.3.2]{DGZ}.
\end{remark}

\begin{theorem}\label{t:ext_lipC1+}
Let $X$ be a super-reflexive Banach space, $E\subset X$, and let $f$ be a real function on $E$.
Consider the following statements:
\begin{enumerate}[(i)]
\item $f$ is bounded and satisfies condition $\cWd G$ for some bounded $G$.
\item $f$ is Lipschitz and satisfies condition $\cWd G$ for some bounded $G$.
\item $f$ can be extended to a Lipschitz function $F\in C^{1,+}(X)$.
\item $f$ can be extended to a function $F\in C^{1,+}(X)$.
\end{enumerate}
Then (i)$\Rightarrow$(ii)$\Leftrightarrow$(iii)$\Rightarrow$(iv).
If $E$ is bounded, then all four statements are equivalent.

Moreover, if (i) or (ii) holds, then $F$ can be found such that $DF(x)=G(x)$ for each $x\in E$.
\end{theorem}
For the most important implication (ii)$\Rightarrow$(iii) the assumption on $X$ cannot be relaxed, see Remark~\ref{r:quasic-bump}.
Cf.~also Remark~\ref{r:ext_loc}.
\begin{proof}
(i)$\Rightarrow$(ii)
Let $K\ge0$ be such that $\norm{G(x)}\le K$ and $\abs{f(x)}\le K$ for each $x\in E$.
Condition $\cWd G$ implies that there is $\de>0$ such that $\absb{big}{f(y)-f(x)-\ev{G(x)}{y-x}}\le\norm{y-x}$ whenever $\norm{y-x}<\de$.
Now let $x,y\in E$.
If $\norm{y-x}<\de$, then
\[
\abs{f(y)-f(x)}\le\absb{big}{f(y)-f(x)-\ev{G(x)}{y-x}}+\absb{big}{\ev{G(x)}{y-x}}\le\norm{y-x}+K\norm{y-x}.
\]
On the other hand, if $\norm{y-x}\ge\de$, then $\abs{f(y)-f(x)}\le2K=\frac{2K}\de\de\le\frac{2K}\de\norm{y-x}$.
Altogether, $f$ is $\max\bigl\{1+K,\frac{2K}\de\bigr\}$-Lipschitz.

(ii)$\Rightarrow$(iii)
By Lemma~\ref{l:WGlip-concave} there exists a bounded concave modulus $\om\in\mc M$ such that $f$ satisfies condition $\cWG\om$ with $G$ on~$E$.
By Theorem~\ref{t:super-norm-modulus} there is an equivalent norm $\tnorm\cdot$ on $X$ that has modulus of smoothness of power type $1+\al$ for some $0<\al\le1$.
Lemma~\ref{l:equiv}(iv) implies that $f$ satisfies $\cWG\om$ in $(X,\tnorm\cdot)$ with the same mapping~$G$.
Denoting $\vp(t)=\int_0^t\om(s)\d s$ for $t\ge0$ and $\nu(x)=\vp(\tnorm x)$ for $x\in X$,
Lemma~\ref{l:smooth_nu}(i) gives that $\nu\in C_{\tnorm\cdot}^{1,\sigma}(X)$ for some bounded $\sigma\in\mc M$.
By Proposition~\ref{p:WGom->extension} used on the space $(X,\tnorm\cdot)$ there is $F\in C_{\tnorm\cdot}^{1,\sigma}(X)$ which is an extension of $f$ with $DF=G$ on $E$.
So $F\in C_{\norm\cdot}^{1,\sigma}(X)\subset C^{1,+}(X)$ by Lemma~\ref{l:equiv}(ii).
Finally, since $\sigma$ is bounded, $DF$ is bounded as well, and therefore $F$ is Lipschitz.

(iii)$\Rightarrow$(ii)
By Fact~\ref{f:conv_subad}, $\om_{DF}\in\mc M$, and so $F$ satisfies condition $\cWG{\om_{DF}}$ with $G_1=DF$ on $X$ by Fact~\ref{f:C1om->WGom}.
Thus $f$ satisfies condition $\cWG{\om_{DF}}$ on $E$ with $G=DF\restr E$ which implies that $f$ satisfies condition $\cWd G$ by~\eqref{e:Wom->WG}.
Moreover, $DF$ is bounded as $F$ is Lipschitz, and so $G$ is bounded and clearly $f=F\restr E$ is Lipschitz.

(iii)$\Rightarrow$(iv) is trivial.

Now suppose that $E$ is bounded and let us prove (iv)$\Rightarrow$(i).
By Fact~\ref{f:conv_subad}, $\om_{DF}\in\mc M$, and so $F$ satisfies condition $\cWG{\om_{DF}}$ on $X$ by Fact~\ref{f:C1om->WGom}.
Thus $f$ satisfies condition $\cWG{\om_{DF}}$ on $E$ with $G=DF\restr E$ which implies that $f$ satisfies condition $\cWd G$ by~\eqref{e:Wom->WG}.
Further, let $B\subset X$ be a ball with $E\subset B$.
Then $DF$ is bounded on $B$, as it is uniformly continuous, and so $F$ is Lipschitz on $B$ and hence bounded on $B$.
Consequently, both $G$ and $f$ are bounded.
\end{proof}

\begin{remark}
In connection with Theorems~\ref{t:ext_C1+_gen} and~\ref{t:ext_lipC1+} note that
condition $\cWd G$ (even with $G=0$) alone does not imply the existence of a $C^{1,+}$-smooth extension even on~$\R$.
Indeed, let $E=\N\subset\R$ and let $f(2n)=0$, $f(2n-1)=n$, and $G(n)=0$ for $n\in\N$.
Then $f$ clearly satisfies $\cWd G$.
However, $f$ cannot be extended to $F\in C^{1,+}(\R)$:
Suppose the contrary.
Then for each $n\in\N$ the Mean value theorem implies the existence of $\zeta_n\in(2n-1,2n)$ and $\xi_n\in(2n,2n+1)$ such that $F'(\zeta_n)=-n$ and $F'(\xi_n)=n+1$.
But $2n+1=F'(\xi_n)-F'(\zeta_n)\le\om_{F'}(2)$ for each $n\in\N$, a contradiction with Fact~\ref{f:conv_subad}.
\end{remark}

\begin{corollary}\label{c:extLipC1+-quasic}
Let $X$ be a super-reflexive Banach space, $c\ge1$, and let $U\subset X$ be an arbitrary union of uniformly separated open $c$-quasiconvex sets.
If $f\in C^{1,+}(U)$ is Lipschitz, then it can be extended to a Lipschitz function $F\in C^{1,+}(X)$.
\end{corollary}
The assumption on $X$ cannot be relaxed even in the most interesting case when we extend from one open bounded quasiconvex set, see Remark~\ref{r:quasic-bump}.
\begin{proof}
Let $\de_0>0$ be such that $U=\bigcup_{\ga\in\Ga}U_\ga$, where each $U_\ga$ is an open $c$-quasiconvex set and $\dist(U_\al,U_\beta)\ge\de_0$ for $\al,\beta\in\Ga$, $\al\neq\beta$.
Denote $f_\ga=f\restr{U_\ga}$.
Given $\ga\in\Ga$, Lemma~\ref{l:quasic-concave_modulus} implies that $\om_{Df_\ga}\in\mc M$ and Lemma~\ref{l:quasic->WG} (with $\om=\om_{Df_\ga}$ and $K=1$)
then implies that $f$ satisfies $\cWG{\om_{Df_\ga}}$ on $U_\ga$ with $G=Df_\ga$ and $M=4c^3$.
We claim that $f$ satisfies condition $\cWd{Df}$.
Indeed, $Df$ is uniformly continuous by the assumption.
Further, let $\ve>0$.
There is $0<\de<\de_0$ such that $\om_{Df}(\de)\le\frac\ve{4c^3}$.
Now if $x,y\in U$, $\norm{y-x}<\de$, then there is $\ga\in\Ga$ such that $x,y\in U_\ga$ and so
$\absb{big}{f(y)-f(x)-\ev{Df(x)}{y-x}}=\absb{big}{f(y)-f(x)-\ev{Df_\ga(x)}{y-x}}\le4c^3\norm{y-x}\om_{Df_\ga}(\norm{y-x})\le4c^3\norm{y-x}\om_{Df}(\norm{y-x})\le\ve\norm{y-x}$.

Moreover, $Df$ is bounded since $f$ is Lipschitz.
Therefore $f$ can be extended to a Lipschitz function $F\in C^{1,+}(X)$ by Theorem~\ref{t:ext_lipC1+}, (ii)$\Rightarrow$(iii).
\end{proof}

\begin{corollary}\label{c:extC1+-quasic}
Let $X$ be a super-reflexive Banach space, $U\subset X$ an open quasiconvex set, and $f\in C^{1,+}(U)$.
Suppose that at least one of the following two conditions is satisfied:
\begin{enumerate}[(a)]
\item $X$ has an equivalent norm with modulus of smoothness of power type~$2$,
\item $U$ is bounded.
\end{enumerate}
Then $f$ can be extended to a function $F\in C^{1,+}(X)$.
\end{corollary}
For the case (b) the assumption of super-reflexivity of $X$ cannot be relaxed, see Remark~\ref{r:quasic-bump}.
\begin{proof}
Suppose that $X$ has an equivalent norm with modulus of smoothness of power type~$2$.
Lemma~\ref{l:quasic-concave_modulus} implies that there exists a concave $\om\in\mc M$ such that $f\in C^{1,\om}(U)$.
Lemma~\ref{l:quasic->WG} then implies that $f$ satisfies $\cWG\om$ on~$U$.
Therefore $f$ can be extended to a function $F\in C^{1,+}(X)$ by Theorem~\ref{t:ext_C1om}.

Now suppose that $U$ is bounded.
Lemma~\ref{l:quasic-concave_modulus} implies that $\om_{Df}\in\mc M$ and so $Df$ is bounded on $U$.
Therefore $f$ is Lipschitz (Proposition~\ref{p:om_der-lip-quasic}) and we may use Corollary~\ref{c:extLipC1+-quasic}.
\end{proof}

\begin{remark}\label{r:ext-lip-quasic}
The assumption of quasiconvexity in Corollaries~\ref{c:extLipC1+-quasic} and \ref{c:extC1+-quasic} (and in Corllary~\ref{c:extC1+B-quasic} below) cannot be dropped even in the case when $U$ is a Jordan domain.
To this end, consider the set
\[
U=U\bigl((0,0),1\bigr)\setminus\bigl\{(x,y)\setsep x\ge0,\abs y\le x^2\bigr\}\subset\R^2
\]
(which coincides with the set $G_2$ from \cite[p.136, Figure 2.1]{BB}) and define $f\colon U\to\R$ by
\[
f(x,y)=\begin{cases}
0\quad&\text{if $x\le0$ or $y<0$,}\\
x^2\quad&\text{if $x\ge0$ and $y\ge0$.}
\end{cases}
\]
Then it is easy to check that $f\in C^{1,\frac12}(U)\subset C^{1,+}(U)$ and that $f$ is Lipschitz.
On the other hand, $f$ cannot be extended to a function $F\in C^1(\R^2)$.
Indeed, suppose that $F$ is a such extension.
Then for all sufficiently small $x>0$ the Mean value theorem implies that for some point $(x,\xi_x)$ with $-x^2\le\xi_x\le x^2$ the equality
\[
\frac{\partial F}{\partial y}(x,\xi_x)=\frac{F(x,x^2)-F(x,-x^2)}{2x^2}=\frac12
\]
holds, which contradicts the fact that $\frac{\partial F}{\partial y}(0,0)=0$.
In particular, $U=G_2$ is not $(C,\om)$-convex for $\om(t)=t^{\frac12}$, contrary to what is stated in~\cite[p.~136]{BB},
since otherwise $f$ would have an extension $F\in C^1(X)$ by \cite[Theorem 2.64]{BB}.
(Cf. also Remark~\ref{r:brud}.)
\end{remark}

\begin{remark}\label{r:quasic_unbound}
We do not know whether conditions (a) and (b) in Corollary~\ref{c:extC1+-quasic} can be disposed of; see Corollary~\ref{c:extC1+B-quasic} for a weaker result.
Note that they can be disposed of if $f$ is Lipschitz (Corollary~\ref{c:extLipC1+-quasic}).
\end{remark}

\begin{remark}\label{r:ext_convex}
We do not know whether it is possible to extend $C^{1,+}$-smooth functions similarly as in Corollaries~\ref{c:extC1om-quasic}, \ref{c:extC1_Hold-quasic}, \ref{c:extLipC1+-quasic}, and \ref{c:extC1+-quasic}
from open bounded \emph{convex} sets (resp. open balls) in some space that is not super-reflexive.
An example of Petr Hájek from \cite{DH} implies that it is not possible in $c_0$ (resp. in a space isomorphic to $c_0$):
There are an open absolutely convex bounded set $U\subset c_0$ and functions $\Phi_\beta\in C^\infty(U)$, $\beta>1$, with all derivatives bounded and Lipschitz,
such that $\Phi_\beta$ cannot be extended to a function from $C^{1,+}(\beta U)$.
This formulation follows rather easily from \cite[remark after Proposition~6.33]{HJ}.
(For a more striking example on a different space see also \cite[Theorem~6.34]{HJ}.)
\end{remark}

\section{\texorpdfstring{The $\Cob{1,+}$ and $\Cloc{1,+}$ cases}{The C_B1+ and C_loc1+ cases}}\label{sec:CB1+}

The following theorem is an extended version of Theorem~\ref{t:ext_CB1+}.

\begin{theorem}\label{t:ext_C1ob}
Let $X$ be a super-reflexive Banach space, $E\subset X$, and let $f$ be a real function on $E$.
Then the following statements are equivalent:
\begin{enumerate}[(i)]
\item $f$ can be extended to a real function $F\in\Cob{1,+}(X)$.
\item $f$ satisfies condition $\cW$ on $E$.
\item For each bounded $A\subset E$ the function $f$ is bounded on $A$ and there exists a bounded $G_A\colon A\to X^*$ such that $f\restr A$ satisfies condition $\cWd{G_A}$ .
\end{enumerate}

Moreover, if $\cW$ is satisfied with $G$, then $F$ can be found such that $DF(x)=G(x)$ for each $x\in E$.
\end{theorem}
For the most important implication (iii)$\Rightarrow$(i) the assumption on $X$ cannot be relaxed, see Remark~\ref{r:quasic-bump}.
\begin{proof}
(i)$\Rightarrow$(ii)
We will show that $f$ satisfies condition $\cW$ on $E$ with $G=DF\restr E$.
To this end consider an arbitrary bounded $B\subset E$ and choose an open convex bounded $U\supset B$.
Denote $H=F\restr U$.
By Fact~\ref{f:conv_subad} and Lemma~\ref{l:Steckin}(iii) there is a concave $\wtilde\om\in\mc M$ such that $H\in C^{1,\wtilde\om}(U)$.
Since $U$ is bounded, $H\in C^{1,\om}(U)$ for $\om(t)=\min\{\wtilde\om(t),\wtilde\om(\diam U)\}$.
As $H$ satisfies $\cWG\om$ with $G_1=DF\restr U$ on $U$ by Fact~\ref{f:C1om->WGom},
Lemma~\ref{l:WGlip-concave} implies that $H$ satisfies $\cWd{G_1}$, $G_1$ is bounded, and $H$ is Lipschitz.
Then $H$ is bounded, since $U$ is bounded.
Consequently, $f\restr B$ satisfies $\cWd{G\restr B}$ and both $G$ and $f$ are bounded on $B$.
Since $B$ was arbitrary, the claim clearly follows.

(ii)$\Rightarrow$(iii) is trivial.

Finally, we will prove (iii)$\Rightarrow$(i) and (ii)$\Rightarrow$(i) with the ``moreover'' part both at once.
If (ii) holds, then we assume that we have a corresponding $G\colon E\to X^*$.
By Theorem~\ref{t:super-norm-modulus} there is an equivalent norm $\nu$ on $X$ with modulus of smoothness of power type $1+\al$ for some $0<\al\le1$.
Then $\nu^{1+\al}\in C_\nu^{1,\al}(X)$ by \cite[proof of Lemma~IV.5.9]{DGZ} combined with Lemma~\ref{l:norm_der_Holder}, and so $\nu^{1+\al}\in C^{1,\al}(X)$ by Lemma~\ref{l:equiv}(ii).
Let $\{\vp_\ga\}_{\ga\in\Ga}$ be a locally finite $C^2$-smooth partition of unity on $\R$ subordinated to the covering $\{(-n,n)\setsep n\in\N\}$
(for the construction of such partition see e.g. \cite[Theorem~3-11]{Spi}).
For each $\ga\in\Ga$ set $\psi_\ga=\vp_\ga\comp\nu^{1+\al}$.
Since the support of each $\vp_\ga$ is compact, it follows that $\vp_\ga'$ is bounded and also $\vp_\ga''$ is bounded, which implies that $\vp_\ga'$ is Lipschitz.
Since $D\nu^{1+\al}$ is $\al$-Hölder on $X$, it follows that $D\nu^{1+\al}$ is bounded on bounded sets and consequently $\nu^{1+\al}$ is Lipschitz on bounded sets (\cite[Proposition~1.71]{HJ}).
Then $\nu^{1+\al}$ is $\al$-Hölder on bounded sets, since for each $r>0$ there is $C>0$ such that $Ct^\al\ge t$ for $t\in[0,r]$.
By \cite[Proposition~1.128(ii)]{HJ} used with $k=1$ and $\om(t)=Kt^\al$ for a sufficiently large $K$ it follows that $\psi_\ga\in C^{1,\al}(U)$ for every bounded open $U\subset X$.

For each $\ga\in\Ga$ the set $\suppo\vp_\ga$ is bounded and so also $\suppo\psi_\ga$ and $\supp\psi_\ga$ are bounded.
Now since for each $\ga\in\Ga$ the function $f\restr{E\cap\supp\psi_\ga}$ is bounded and satisfies condition $\cWd H$ for some bounded $H\colon E\cap\supp\psi_\ga\to X^*$ (by (iii) with $A=E\cap\supp\psi_\ga$),
Theorem~\ref{t:ext_lipC1+} implies that there exists $F_\ga\in C^{1,+}(X)$ which is an extension of $f\restr{E\cap\supp\psi_\ga}$.
In case that (ii) holds we can take $H=G\restr{E\cap\supp\psi_\ga}$ and so by Theorem~\ref{t:ext_lipC1+} we may additionally assume that $DF_\ga(x)=G(x)$ for each $x\in E\cap\supp\psi_\ga$.
Set $F=\sum_{\ga\in\Ga}\psi_\ga F_\ga$ and note that the sum is in fact finite on bounded sets:
Let $r>0$.
Since $\{\suppo\vp_\ga\}$ is a locally-finite system, it follows that there is $\Ga_r\subset\Ga$ finite such that $\vp_\ga=0$ on the compact interval $[0,r^{1+\al}]$ for $\ga\notin\Ga_r$.
Hence $\psi_\ga=0$ (and so also $\psi_\ga F_\ga=0$) on $U_\nu(0,r)$ for $\ga\notin\Ga_r$, which implies that $F(x)=\sum_{\ga\in\Ga_r}\psi_\ga(x)F_\ga(x)$ for $x\in U_\nu(0,r)$.

Further, note that $\psi_\ga F_\ga\in\Cob{1,+}(X)$ for each $\ga\in\Ga$.
Indeed, if $V\subset X$ is an open ball, then $\psi_\ga$, $D\psi_\ga$, and $DF_\ga$ are bounded on~$V$.
It follows that $\psi_\ga$ and $F_\ga$ are Lipschitz on $V$, thus $F_\ga$ is bounded on $V$,
and in turn all four mappings $\psi_\ga$, $D\psi_\ga$, $F_\ga$, and $DF_\ga$ are uniformly continuous on $V$ with the same modulus (Fact~\ref{f:common_mod}).
It follows that $\psi_\ga F_\ga\in C^{1,+}(V)$ (\cite[Proposition~1.129]{HJ}).
Consequently, $F\in\Cob{1,+}(X)$, as the sum in its definition is finite on bounded sets.

To see that $F$ is an extension of $f$ let $x\in E$.
Let $r>0$ be such that $x\in U_\nu(0,r)$ and let $\Ga_r$ be as above.
Put $\Lambda=\{\ga\in\Ga_r\setsep x\in\supp\psi_\ga\}$.
Then $\psi_\ga(x)=0$ for $\ga\in\Ga\setminus\Lambda$ and $F_\ga(x)=f\restr{E\cap\supp\psi_\ga}(x)=f(x)$ whenever $\ga\in\Lambda$.
Hence
\[
F(x)=\sum_{\ga\in\Lambda}\psi_\ga(x)F_\ga(x)=\sum_{\ga\in\Lambda}\psi_\ga(x)f(x)=f(x)\sum_{\ga\in\Lambda}\psi_\ga(x)=f(x)\sum_{\ga\in\Ga}\psi_\ga(x)=f(x).
\]

Finally, assuming that (ii) holds we show that $DF(x)=G(x)$.
Since $\sum_{\ga\in\Ga_r}\psi_\ga(y)=1$ for any $y\in U_\nu(0,r)$, it follows that $\sum_{\ga\in\Ga}D\psi_\ga(x)=\sum_{\ga\in\Ga_r}D\psi_\ga(x)=0$.
Moreover, $D\psi_\ga(x)=0$ whenever $\ga\in\Ga_r\setminus\Lambda$, as the set $X\setminus\supp\psi_\ga$ is open,
and $DF_\ga(x)=G(x)$ for $\ga\in\Lambda$.
Hence
\[\begin{split}
DF(x)&=\sum_{\ga\in\Ga_r}\psi_\ga(x)DF_\ga(x)+F_\ga(x)D\psi_\ga(x)=\sum_{\ga\in\Lambda}\psi_\ga(x)DF_\ga(x)+F_\ga(x)D\psi_\ga(x)\\
&=\sum_{\ga\in\Lambda}\psi_\ga(x)G(x)+f(x)D\psi_\ga(x)=G(x)\sum_{\ga\in\Ga}\psi_\ga(x)+f(x)\sum_{\ga\in\Ga}D\psi_\ga(x)=G(x).
\end{split}\]
\end{proof}

\begin{corollary}\label{c:extC1+B-quasic}
Let $X$ be a super-reflexive Banach space, $U\subset X$ an open quasiconvex set, and $f\in C^{1,+}(U)$.
Then $f$ can be extended to a function $F\in\Cob{1,+}(X)$.
\end{corollary}
The assumption on $X$ cannot be relaxed, see Remark~\ref{r:quasic-bump}.
\begin{proof}
Lemma~\ref{l:quasic-concave_modulus} implies that there exists $\om\in\mc M$ such that $f\in C^{1,\om}(U)$.
Lemma~\ref{l:quasic->WG} then implies that $f$ satisfies $\cWG\om$ on $U$ with $G=Df$ and some $M>0$.
Now let $A\subset U$ be an arbitrary bounded subset.
Fixing some $x_0\in A$, the inequalities in $\cWG\om$ imply that $\norm{G(y)}\le\norm{G(x_0)}+M\om(\diam A)$ and $\abs{f(y)}\le\abs{f(x_0)}+\norm{G(x_0)}\diam A+M\om(\diam A)\diam A$ for each $y\in A$,
and so both $f$ and $G$ are bounded on $A$.
Also, $f\restr A$ satisfies $\cWd{G\restr A}$ by~\eqref{e:Wom->WG}.
Therefore $f$ can be extended to a function $F\in\Cob{1,+}(X)$ by Theorem~\ref{t:ext_C1ob}.
\end{proof}

At least for certain quasiconvex sets it is possible to relax the assumption on $f$ in the previous corollary:
\begin{corollary}\label{c:extC1+B-convex}
Let $X$ be a super-reflexive Banach space, let $U\subset X$ be an open set that is a bi-Lipschitz image of a convex subset of a normed linear space, and let $f\in\Cob{1,+}(U)$.
Then $f$ can be extended to a function $F\in\Cob{1,+}(X)$.
\end{corollary}
\begin{proof}
Let $\Phi\colon U\to V$ be a bi-Lipschitz mapping onto, where $V$ is a convex subset of some normed linear space.
Let $A\subset U$ be an arbitrary bounded subset.
Then $\Phi(A)$ is bounded and so there is $R>0$ such that $\Phi(A)\subset U(0,R)$.
The set $V\cap U(0,R)$ is convex and relatively open in $V$, and so $W=\Phi^{-1}(V\cap U(0,R))$ is open and it is a bi-Lipschitz image of a bounded convex set, hence it is quasiconvex (Remark~\ref{r:biLip-quasic}) and bounded.
Moreover, $A\subset W$.

Since $f\restr W\in C^{1,+}(W)$, Corollary~\ref{c:extC1+-quasic} combined with Theorem~\ref{t:ext_lipC1+} implies that $f\restr W$ is bounded and satisfies condition $\cWd G$ for some bounded~$G$.
Thus $f$ is bounded on $A$ and $f\restr A$ satisfies condition $\cWd{G\restr A}$ for some bounded~$G$.
Therefore $f$ can be extended to a function $F\in\Cob{1,+}(X)$ by Theorem~\ref{t:ext_C1ob}.
\end{proof}

\begin{remark}\label{r:quasic-bump}
In Remark~\ref{r:part_unity} we showed that the assumptions on $X$ in Theorems~\ref{t:ext_C1om} and \ref{t:extC1_Holder} cannot be relaxed.
Now we show that
\begin{equation}\label{e:assum}\begin{split}
&\text{the assumption on $X$ in Corollaries~\ref{c:extC1om-quasic}, \ref{c:extC1_Hold-quasic}, \ref{c:extLipC1+-quasic}, \ref{c:extC1+-quasic}``(b)'', \ref{c:extC1+B-quasic}, and Theorems~\ref{t:ext_lipC1+} and \ref{t:ext_C1ob} cannot be relaxed.}
\end{split}\end{equation}

(a)
We start with a construction of a bounded open quasiconvex set $U$ in an arbitrary normed linear space $X$ with $\dim X\ge2$ and a $C^{1,1}$-smooth function $f$ on it.
Fix $e\in S_X$.
Set $C=\bigcup_{t>0}tU(e,\frac14)$ and $U=(U(0,2)\setminus B_X)\cup(U(0,2)\cap C)$.
We claim that $U$ is $5$-quasiconvex.
Since $U(0,2)\setminus B_X$ is $2$-quasiconvex by Example~\ref{ex:shell-quasic} and $C$ is convex, it suffices to consider $x\in U(0,2)\setminus B_X$ and $y\in C\cap B_X$.
Set $z=\norm x\frac y{\norm y}\in U(0,2)\setminus B_X$.
Then $\norm{z-y}=\norm z-\norm y=\norm x-\norm y\le\norm{x-y}$ and $\norm{x-z}\le\norm{x-y}+\norm{y-z}\le2\norm{x-y}$.
By the $2$-quasiconvexity of $U(0,2)\setminus B_X$ there is a continuous rectifiable curve $\ga$ joining $x$ and $z$ in $U$ such that $\len\ga\le2\norm{x-z}$.
By concatenating $\ga$ with the segment $[z,y]\subset U$ we obtain a continuous rectifiable curve joining $x$ and $y$ in $U$ whose length is at most $2\norm{x-z}+\norm{z-y}\le5\norm{x-y}$.

Next, by the Hahn-Banach theorem there is $g\in S_{X^*}$ such that $g(e)=1$.
Note that
\begin{equation}\label{e:norma_x}
\text{$\norm x<\frac56$ whenever $x\in C$ and $g(x)\le\frac12$.}
\end{equation}
Indeed, there is $t>0$ such that $x\in U(te,\frac14t)$.
Then $t=g(te-x)+g(x)\le\norm{te-x}+\frac12<\frac14t+\frac12$, which implies that $t<\frac23$, and consequently $\norm x\le\norm{x-te}+\norm{te}<\frac54t<\frac56$.

Now let $\vp\in C^{1,1}(\R)$ be such that $\vp=0$ on $[\frac12,+\infty)$, $\vp(\frac14)>0$, and $\abs{\vp'}$ is bounded by a constant $M>0$.
Define $f\colon U\to\R$ by $f=\vp\comp g$ on $C\cap U$ and $f=0$ on $U\setminus B_X$.
Notice that $f$ is well-defined, since $g>\frac12$ on $C\setminus B_X$ by~\eqref{e:norma_x}.
The function $f$ is clearly $C^{1,1}$-smooth on $C\cap U$ and also on $U\setminus\{x\in C\setsep g(x)\le\frac12\}$ ($Df=0$ on this set).
Further, if $x\in C\cap U$ with $g(x)\le\frac12$ and $y\in U\setminus C$, then $\norm{Df(x)-Df(y)}=\norm{Df(x)}=\norm{D\vp(g(x))\comp g}\le\abs{\vp'(g(x))}\cdot\norm g\le M$.
On the other hand, $\norm{x-y}\ge\norm y-\norm x>1-\frac56=\frac16$ (use~\eqref{e:norma_x}).
It follows that $\norm{Df(x)-Df(y)}\le M=6M\frac16<6M\norm{x-y}$.
Consequently, $f\in C^{1,1}(U)$.

Note that $f\in C^{1,1}(U)\subset C^{1,\al}(U)\subset C^{1,+}(U)$ for any $0<\al\le1$ (the first inclusion follows since $U$ is bounded) and that $f$ is Lipschitz (Proposition~\ref{p:om_der-lip-quasic}).
Let $S$ be one of the spaces $C^{1,\al}(X)$, $C^{1,+}(X)$, $\Cob{1,+}(X)$.
Assume that we can extend $f$ to a function $F\in S$.
Let $H=\chi_{U(0,2)}\cdot F$.
Since $H=0$ on $X\setminus B_X$ and also $H=F=0$ on $U(0,2)\setminus B_X$, it is easy to check that $H\in S$.
Moreover, $H$ is a bump function.
In particular, $H\in C^{1,+}(X)$ even in the case $S=\Cob{1,+}(X)$.
It follows that $X$ is super-reflexive and in case that $S=C^{1,\al}(X)$ it has an equivalent norm with modulus of smoothness of power type $1+\al$ (\cite[Theorem~V.3.2]{DGZ}).

\smallskip
(b)
To prove~\eqref{e:assum} assume that Corollary~\ref{c:extC1om-quasic}, resp. \ref{c:extC1_Hold-quasic}, resp. \ref{c:extLipC1+-quasic}, resp. \ref{c:extC1+-quasic}(b), resp.~\ref{c:extC1+B-quasic}
holds for some $X$ which does not satisfy the respective assumptions.
Now we apply (a) with $S=C^{1,1}(X)=C^{1,\om}(X)$ for $\om(t)=t$ for Corollary~\ref{c:extC1om-quasic}, resp. $S=C^{1,\al}(X)$ for Corollary~\ref{c:extC1_Hold-quasic},
resp. $S=C^{1,+}(X)$ for Corollaries~\ref{c:extLipC1+-quasic}, \ref{c:extC1+-quasic}, resp. $S=\Cob{1,+}(X)$ for Corollary~\ref{c:extC1+B-quasic}, which leads to a contradiction.
Further, since the assumptions on $X$ in Corollary~\ref{c:extLipC1+-quasic}, resp. \ref{c:extC1+B-quasic}, cannot be relaxed,
it follows that the assumptions on $X$ in Theorem~\ref{t:ext_lipC1+}, resp. \ref{t:ext_C1ob}, cannot be relaxed as well,
since Corollary~\ref{c:extLipC1+-quasic}, resp. \ref{c:extC1+B-quasic}, is a direct consequence of Theorem~\ref{t:ext_lipC1+}, resp.~\ref{t:ext_C1ob}.
\end{remark}

Using Theorem~\ref{t:ext_lipC1+} and another partition of unity we obtain the following result.
\begin{theorem}\label{t:ext_C1+loc}
Let $X$ be a super-reflexive Banach space, $U\subset X$ an open set, $A\subset U$ relatively closed in $U$, and $f\colon A\to\R$.
Then $f$ can be extended to a function $F\in\Cloc{1,+}(U)$ if and only if
$f$ is locally bounded and for each $x\in A$ there exist $\de_x>0$ and a bounded $G_x\colon A\cap U(x,\de_x)\to X^*$ such that $f\restr{A\cap U(x,\de_x)}$ satisfies condition $\cWd{G_x}$.
\end{theorem}
\begin{proof}
$\Rightarrow$
$F$ is clearly locally bounded and hence so is $f$.
Choose any $x\in A$.
There is $\de_x>0$ such that $DF$ is uniformly continuous on $U(x,\de_x)\subset U$.
Set $H=F\restr{U(x,\de_x)}$.
Since $\om_{DH}\in\mc M$ by Fact~\ref{f:conv_subad}, $DH$ is bounded and Fact~\ref{f:C1om->WGom} implies that $H$ satisfies condition $\cWG{\om_{DH}}$ on $U(x,\de_x)$ with $G=DH$.
Therefore $H$ clearly satisfies condition $\cWd{DH}$ (see~\eqref{e:Wom->WG}).
Consequently, $f\restr{A\cap U(x,\de_x)}$ satisfies condition $\cWd{DH\restr{A\cap U(x,\de_x)}}$.

$\Leftarrow$
We may assume that $U(x,\de_x)\subset U$ and $f$ is bounded on $A\cap U(x,\de_x)$ for every $x\in A$.
Let $\{\vp_\al\}_{\al\in\Lambda'}$ be a locally finite $C^{1,+}$-smooth partition of unity on $U$ subordinated to the open covering $\{U(x,\de_x)\}_{x\in A}\cup\{U\setminus A\}$.
(The existence of such partition follows from Theorem~\ref{t:super-norm-modulus} and \cite[Proposition~3.7]{Kr2}, whose proof is a simple modification of the proof of~\cite[Theorem~7.56]{HJ}.
Another possibility is described in Remark~\ref{r:part_unity}.)
Set $\Lambda=\{\al\in\Lambda'\setsep\suppo\vp_\al\cap A\neq\emptyset\}$.
For each $\al\in\Lambda$ let $x_\al\in A$ be such that $\suppo\vp_\al\subset U(x_\al,\de_{x_\al})$.
Denote $U_\al=U(x_\al,\de_{x_\al})$.
By Theorem~\ref{t:ext_lipC1+} for each $\al\in\Lambda$ there is $F_\al\in C^{1,+}(X)$ which is an extension of $f\restr{A\cap U_\al}$.
Now put $F=\sum_{\al\in\Lambda}\vp_\al F_\al$.

First note that $\vp_\al F_\al\in\Cloc{1,+}(U)$ for each $\al\in\Lambda$.
Indeed, given $x\in U$, by continuity there is $\de>0$ such that $V=U(x,\de)\subset U$ and $\vp_\al$, $D\vp_\al$, $F_\al$, and $DF_\al$ are bounded on $V$.
It follows that $\vp_\al$ and $F_\al$ are Lipschitz on $V$, and in turn all four mappings are uniformly continuous on $V$ with the same modulus (Fact~\ref{f:common_mod}).
It follows that $\vp_\al F_\al\in C^{1,+}(V)$ (\cite[Proposition~1.129]{HJ}).
Consequently, $F\in\Cloc{1,+}(U)$, as the sum in its definition is locally finite:
for every $x\in U$ there is a neighbourhood $W$ of $x$ and $H\subset\Lambda$ finite such that $\vp_\al F_\al=0$ on $W$ for $\al\in\Lambda\setminus H$, and so $F=\sum_{\al\in H}\vp_\al F_\al$ on~$W$.

Finally, to show that $F$ is an extension of $f$ suppose that $x\in A$ is given.
Then $\vp_\al(x)=0$ for each $\al\in\Lambda'\setminus\Lambda$ and for each $\al\in\Lambda$ such that $x\notin U_\al$.
Hence
\[
F(x)=\sum_{\substack{\al\in\Lambda\\x\in U_\al}}\vp_\al(x)F_\al(x)=\sum_{\substack{\al\in\Lambda\\x\in U_\al}}\vp_\al(x)f(x)=f(x)\sum_{\al\in\Lambda'}\vp_\al(x)=f(x).
\]
\end{proof}

\begin{remark}\label{r:ext_loc}
The assumption on $X$ in the previous theorem cannot be relaxed even in the case $U=X$.
In fact we show more:
Even a weaker version of (i)$\Rightarrow$(iv) in Theorem~\ref{t:ext_lipC1+}, namely
\begin{equation}\label{e:weaker}
\text{(i) implies that $f$ can be extended to $F\in\Cloc{1,+}(X)$}
\end{equation}
implies that $X$ is super-reflexive.
Note that the validity of Theorem~\ref{t:ext_C1+loc} for some space $X$ implies~\eqref{e:weaker},
since the ``global'' statement~(i) in Theorem~\ref{t:ext_lipC1+} is stronger than the condition in Theorem~\ref{t:ext_C1+loc}, which is a ``localised'' version of this statement.

Indeed, suppose that~\eqref{e:weaker} holds for some space $X$.
Then there exists a $\Cloc{1,+}$-smooth bump on $X$ (extend the function defined by $f=0$ on $X\setminus U(0,1)$ and $f(0)=1$).
We claim that $X$ does not contain a subspace isomorphic to $c_0$ and so $X$ is super-reflexive by \cite{FWZ}, see \cite[Corollary~5.51]{HJ}.
Suppose to the contrary that $X$ contains a subspace $Y$ isomorphic to $c_0$.
By the construction below there is a bounded function $f$ on a closed subset $A\subset c_0$ that satisfies $\cWd G$ with $G=0$, where $G\colon A\to c_0^*$, and yet
\begin{equation}\label{e:yet}
\text{$f$ cannot be extended to a function on $c_0$ that is $C^{1,+}$-smooth on a neighbourhood of $0$.}
\end{equation}
Identifying $A$ with the corresponding subset of $Y$ it is easy to see that $f$ satisfies condition $\cWd H$ in $X$ with $H=0$, where $H\colon A\subset Y\to X^*$.
So by~\eqref{e:weaker} the function $f$ can be extended to $F\in\Cloc{1,+}(X)$, which contradicts~\eqref{e:yet}.

We set $A=\{0\}\cup\bigl\{\frac1{2^n}e_k\setsep n,k\in\N\bigr\}$ and define $f\colon A\to\R$ by $f(0)=0$ and $f\bigl(\frac1{2^n}e_k\bigr)=\bigl(\frac1{2^n}\bigl)^2$.
The set $A$ is clearly closed.
We claim that $f$ satisfies $\cWd G$ with $G=0$:
Given $\ve>0$ set $\de=\frac\ve4$.
Now let $x,y\in A$, $0<\norm{x-y}<\de$.
We may assume that $x=\frac1{2^n}e_k$ for some $n,k\in\N$.
If $y=0$, then $\abs{f(x)-f(y)}=\bigl(\frac1{2^n}\bigl)^2=\norm{x-y}^2<\de\norm{x-y}<\ve\norm{x-y}$.
Otherwise $y=\frac1{2^m}e_l$ for some $m,l\in\N$.
We may assume without loss of generality that $m\ge n$.
If $k=l$, then $\norm{x-y}=\frac1{2^n}-\frac1{2^m}\ge\frac1{2^{n+1}}$, otherwise $\norm{x-y}=\frac1{2^n}>\frac1{2^{n+1}}$.
Therefore $\abs{f(x)-f(y)}<\bigl(\frac1{2^n}\bigl)^2\le4\norm{x-y}^2<4\de\norm{x-y}=\ve\norm{x-y}$.

Now since $\lim_{k\to\infty}\frac1{2^n}e_k=0$ weakly, it follows that $f$ is not weakly sequentially continuous on any neighbourhood of~$0$.
Thus \eqref{e:yet} follows from results in \cite{Haj:Smooth-c0}, see \cite[Theorem~6.30]{HJ} and notice that clearly $\mathcal C_{\mathrm{wsC}}\subset\mathcal C_{\mathrm{wsc}}$ (cf. \cite[pp.~137--138]{HJ}).
\end{remark}

We would like to thank Petr Hájek for providing us with the example in Remark~\ref{r:ext_loc}.

\end{document}